\documentclass[reqno, 12pt]{amsart}

\usepackage{amsmath}
\usepackage{amsfonts}
\usepackage{amssymb}
\usepackage{mathrsfs}
\usepackage{amsthm}
\usepackage{bbm}
\usepackage{graphicx, color}
\usepackage{graphbox}
\usepackage{tikz}
\usetikzlibrary{cd, intersections, calc, decorations.pathmorphing, arrows, decorations.pathreplacing}
\usepackage{array}
\DeclareMathAlphabet{\mathpzc}{OT1}{pzc}{m}{it}

\usepackage{setspace}
\usepackage{caption}
\usepackage{here}

\usepackage{mathtools}

\usepackage{enumerate}

\definecolor{refkey}{rgb}{0.9451,0.2706,0.4941}
\definecolor{labelkey}{rgb}{0.9451,0.2706,0.4941}

\definecolor{mygreen}{rgb}{0,0.7,0.3}
\definecolor{myblue}{rgb}{0,0.50,1.20}
\definecolor{myorange}{rgb}{1,0.38,0}

\usepackage{subfiles}
\usepackage[colorlinks, linkcolor=blue, citecolor=red, urlcolor=cyan,pagebackref]{hyperref}

\numberwithin{equation}{section}

\allowdisplaybreaks

\usepackage{cleveref}
\crefname{thm}{Theorem}{Theorems}
\crefname{introthm}{Theorem}{Theorems}
\crefname{cor}{Corollary}{Corollaries}
\crefname{lem}{Lemma}{Lemmas}
\crefname{sublem}{Sublemma}{Sublemmas}
\crefname{prop}{Proposition}{Propositions}
\crefname{dfn}{Definition}{Definitions}
\crefname{defi}{Definition}{Definitions}
\crefname{ex}{Example}{Examples}
\crefname{claim}{Claim}{Claims}
\crefname{conj}{Conjecture}{Conjectures}
\crefname{conv}{Notation}{Notations}
\crefname{rem}{Remark}{Remarks}
\crefname{rmk}{Remark}{Remarks}
\crefname{prob}{Problem}{Problems}
\crefname{figure}{Figure}{Figures}
\crefname{table}{Table}{Tables}
\crefname{section}{Section}{Sections}
\crefname{subsection}{Section}{Sections}
\crefname{appendix}{Appendix}{Appendices}

\newtheorem{thm}{Theorem}[section]
\newtheorem{prop}[thm]{Proposition}
\newtheorem{cor}[thm]{Corollary}
\newtheorem{lem}[thm]{Lemma}

\theoremstyle{definition}
\newtheorem{dfn}[thm]{Definition}
\newtheorem{defi}[thm]{Definition}

\theoremstyle{remark}
\newtheorem{rmk}[thm]{Remark}
\newtheorem{rem}[thm]{Remark}
\newtheorem{prob}[thm]{Problem}

\newcommand*{\chom}{\mathcal{H}\kern -.5pt om}

\newcommand{\bZ}{\mathbb{Z}}
\newcommand{\bQ}{\mathbb{Q}}
\newcommand{\bR}{\mathbb{R}}

\newcommand{\bS}{\mathbb{S}}

\newcommand{\bP}{\mathbb{P}}

\newcommand{\bE}{\mathbb{E}}

\newcommand{\bs}{{\boldsymbol{s}}}

\newcommand{\HHH}{\mathbb{H}}
\newcommand{\SSS}{\mathbb{S}_{\infty}^1}

\newcommand{\cA}{\mathcal{A}}

\newcommand{\cC}{\mathcal{C}}
\newcommand{\cD}{\mathcal{D}}
\newcommand{\cE}{\mathcal{E}}
\newcommand{\cF}{\mathcal{F}}

\newcommand{\cQ}{\mathcal{Q}}

\newcommand{\cS}{\mathcal{S}}
\newcommand{\cT}{\mathcal{T}}

\newcommand{\X}{\mathcal{X}}
\newcommand{\cX}{\mathcal{X}}

\newcommand{\fF}{\mathfrak{F}}

\newcommand{\fS}{\mathfrak{S}}

\newcommand{\tri}{\triangle}
\newcommand{\sgn}{\mathrm{sgn}}
\newcommand{\trop}{\mathrm{trop}}
\newcommand{\pos}{\mathbb{R}_{>0}}

\newcommand{\ML}{\mathcal{ML}}

\newcommand{\eML}{\widehat{\mathcal{ML}}}

\newcommand{\tr}{\mathsf{T}}

\newcommand{\bExch}{\bE \mathrm{xch}}

\newcommand{\Teich}{Teichm\"uller}
\newcommand{\cl}{\colon}
\newcommand{\ceq}{\coloneqq}
\newcommand{\PSL}{PSL(2,\bR)}

\newcommand{\dsp}{\displaystyle}
\newcommand{\ve}{\varepsilon}

\newcommand{\fin}{\mathrm{pre}}

\DeclareMathOperator{\interior}{\mathrm{int}}

\DeclareMathOperator{\Frac}{\mathrm{Frac}}

\makeatletter
\newcommand{\oset}[3][0ex]{%
  \mathrel{\mathop{#3}\limits^{
    \vbox to#1{\kern-2\ex@
    \hbox{$\scriptstyle#2$}\vss}}}}
\makeatother
\newcommand{\overbar}[1]{\oset{#1}{-\!\!\!-\!\!\!-}}

\makeatletter
\newcommand{\osetnear}[3][0ex]{%
  \mathrel{\mathop{#3}\limits^{
    \vbox to#1{\kern-.3\ex@
    \hbox{$\scriptstyle#2$}\vss}}}}
\makeatother
\newcommand{\overbarnear}[1]{\osetnear{#1}{-\!\!\!-\!\!\!-}}

\newcommand\qarrow[2]{\draw[->,shorten >=4pt,shorten <=4pt] (#1) -- (#2) [thick];} 

\tikzset{
  mid arrow/.style={postaction={decorate,decoration={
        markings,
        mark=at position .5 with {\arrow[#1]{stealth}}
      }}},
}

\tikzset{
  symbol/.style={
    draw=none,
    every to/.append style={
      edge node={node [sloped, allow upside down, auto=false]{$#1$}}}
  }
}

\renewcommand{\mathbf}{\boldsymbol}

\title{Earthquake theorem for cluster algebras \protect\linebreak of finite type}

\author[Takeru Asaka]{Takeru Asaka}
\address{Takeru Asaka, Graduate School of Information Sciences, Tohoku University,  6-3 Aoba, Aramaki, Aoba-ku, Sendai, Miyagi 980-8578, Japan.}
\email{takeru.asaka.b5@tohoku.ac.jp}

\author[Tsukasa Ishibashi]{Tsukasa Ishibashi}
\address{Tsukasa Ishibashi, Mathematical Institute, Tohoku University, 6-3 Aoba, Aramaki, Aoba-ku, Sendai, Miyagi 980-8578, Japan.}
\email{tsukasa.ishibashi.a6@tohoku.ac.jp}

\author[Shunsuke Kano]{Shunsuke Kano}
\address{Shunsuke Kano, Mathematical Science Center for Co-creative Society, Tohoku University
6-3 Aoba, Aramaki, Aoba-ku, Sendai, Miyagi 980-8578, Japan.}
\email{s.kano@tohoku.ac.jp}

\date{\today}

\begin{document}
\maketitle

\begin{abstract}
We introduce a cluster algebraic generalization of Thurston's earthquake map for the cluster algebras of finite type, which we call the \emph{cluster earthquake map}. 
It is defined by gluing exponential maps, which is modeled after the earthquakes along ideal arcs. 
We prove an analogue of the earthquake theorem, which states that the cluster earthquake map gives a homeomorphism between the spaces of $\mathbb{R}^\mathrm{trop}$- and $\mathbb{R}_{>0}$-valued points of the cluster $\mathcal{X}$-variety. 
For those of type $A_n$ and $D_n$, 
the cluster earthquake map indeed recovers the 
earthquake maps for marked disks and once-punctured marked disks, respectively. 
Moreover, we investigate certain asymptotic behaviors of the cluster earthquake map, which give rise to ``continuous deformations'' of the Fock--Goncharov fan.
\end{abstract}

\tableofcontents

\section{Introduction}

\subsection{Thurston's earthquake maps}

Earthquakes are particular deformations of hyperbolic structures on a closed surface $\Sigma$ along measured geodesics laminations, which is introduced by W.~P.~Thurston \cite{Thu84}. 
For example, the (left) earthquake along a simple closed curve $C$ with a transverse measure $\mu$ deforms a given hyperbolic structure by cutting the surface along $C$ and gluing back after rotating the ``right" side of $C$ by $\mu$. 
Bonsante--Krasnov--Shlenker \cite{BKS16} generalized the earthquakes to the setting of a marked surface, which is a compact oriented surface $\Sigma$ with a fixed finite set $M \subset \Sigma$ of marked points (see \cref{cl_surf}).
Earthquakes give a map $E: \widehat{\cT}(\Sigma) \times \eML(\Sigma) \to \widehat{\cT}(\Sigma)$, where $\widehat{\cT}(\Sigma)$ is an \emph{enhanced \Teich\ space} and $\eML(\Sigma)$ is a certain space of measured geodesic laminations: see \cref{cl_surf}. It is equivariant under the mapping class group.
Most notably, the \emph{earthquake theorem} states that any two hyperbolic structures are related by a unique earthquake. In other words, for any point $g_0 \in \widehat{\cT}(\Sigma)$, the map 
\begin{align}\label{eq:Eq_homeo}
    E(g_0,-): \eML(\Sigma) \xrightarrow{\sim} \widehat{\cT}(\Sigma)
\end{align}
gives a homeomorphism.
This is proved by \cite{Thu84} for a closed surface, and by \cite[Theorem 1.4]{BKS16} for a closed surface with punctures. 
Among its applications, Kerckhoff \cite{Ker83} succeeded in solving the Nielsen realization problem of finite subgroups of mapping class groups. 
The map \eqref{eq:Eq_homeo} can be viewed as a kind of ``exponential map'' based at $g_0$. Indeed, the coordinate description of the earthquakes along ideal arcs take the form of exponential map. See \cref{eq_surf}.

After a couple of decades, Fock--Goncharov \cite{FG07} pointed out that the space $\eML(\Sigma)$ can be regarded as a ``tropical analogue'' of $\widehat{\cT}(\Sigma)$. Their observation is that the coordinate transformation formula between the shear coordinates on $\eML(\Sigma)$ is the tropical limit of those between the cross ratio coordinates on $\widehat{\cT}(\Sigma)$ when we change ideal triangulations of $\Sigma$. 
In the language of \emph{cluster varieties} \cite{FG09}, it can be rephrased that these spaces are canonically isomorphic to the spaces of real tropical points $\bR^\trop$ and positive real points $\pos$, respectively, of a cluster variety $\X_{\Sigma}$ arising from the moduli space of framed $PGL_2$-local systems on $\Sigma$. Symbolically, 
\begin{align}\label{eq:ML_T_cluster}
    \eML(\Sigma) \cong \X_\Sigma(\bR^\trop)\quad \mbox{and} \quad \widehat{\cT}(\Sigma)\cong \X_\Sigma(\pos).
\end{align}
The theory of cluster varieties is the geometric counterpart of the theory of \emph{cluster algebras} initiated by Fomin--Zelevinsky \cite{FZ-CA1}. 
Since its appearance, the theory of cluster varieties/algebras has been rapidly developed, tightly related to the representation theory of algebraic groups \cite{FZ-CA2,BFZ,Williams} finite-dimensional algebras and quivers \cite{BMRRT,CC}, higher \Teich\ theory \cite{FG03,GS19}, integrable systems \cite{GK13,FH14}, and so on. These references are far from exhaustive. A good overview on the algebraic aspect is found in \cite{Nak21-1}.  
\subsection{Cluster earthquake maps}
Our aim in this paper is to explore the cluster algebraic aspects of the earthquake theorem. In particular, we ask if a cluster algebraic analogue of the earthquake theorem can be formulated and proved. 
Firstly, in view of the isomorphisms \eqref{eq:ML_T_cluster}, the earthquake map $E$ can be rewritten as 
\begin{align*}
    E: \cX_\Sigma (\pos)\times\cX_\Sigma(\bR^\trop)\to \cX_\Sigma (\pos),
\end{align*}
and the earthquake theorem is rephrased that the map
\begin{align*}
    E(g_0,-): \X_\Sigma(\bR^\trop) \xrightarrow{\sim} \X_\Sigma(\pos)
\end{align*}
is a homeomorphism for any point $g_0 \in \X_\Sigma(\pos)$. Moreover, based on \cite{A18}, we see that the restriction $E(g_0,-)|_{\cC^+_\tri}$ takes the form of the exponential map (\cref{eq_surf}), where
$\cC^+_{\tri} \subset \X_\Sigma(\bR^\trop)$ is the cone consisting of the laminations with non-negative shear coordinates associated with an ideal triangulation $\tri$.
Hence, the restriction 
\begin{align*}
    E_\fin(g_0,-): |\fF^+_\Sigma| \to \X_\Sigma(\pos)
\end{align*}
of $E(g_0,-)$ to the support $|\fF^+_\Sigma|=\bigcup_{\tri} \cC^+_\tri \subset \X_\Sigma(\bR^\trop)$ of the Fock--Goncharov fan (\cref{def:positive_cone})
is something that might be called a ``piecewise exponential map''. Let us call the map $E_\fin: |\fF^+_\Sigma| \times \X_\Sigma(\pos) \to \X_\Sigma(\pos)$ the \emph{pre-earthquake map}. 
In the surface case, it is known that $|\fF^+_\Sigma|$ is dense in $\X_\Sigma(\bR^\trop)$ \cite{Yurikusa} and hence the earthquake map can be viewed as a continuous extension of the pre-earthquake map. 

Generalizing the above observation, we define the \emph{cluster pre-earthquake map} (\cref{def:exp})
\begin{align}\label{eq:pre-Eq}
    E_\fin: \X_\bs(\pos)\times  |\fF^+_\bs| \to \X_\bs(\pos)
\end{align}
for the cluster variety $\X_\bs$ associated with any mutation class $\bs$ of seeds. Here $\fF^+_\bs$ denotes the Fock--Goncharov fan. The map $E_\fin$ is equivariant under the action of the cluster modular group $\Gamma_\bs$ (\cref{lem:equivariance}). 
Let us consider the mutation class $\bs$ of finite type (namely, there are only finitely many seeds in $\bs$), which are known to be classified by Dynkin diagrams \cite{FZ-CA2}. 
In the finite type case, it is known that the fan $\fF^+_\bs$ is complete: $|\fF^+_\bs|=\X_\bs(\bR^\trop)$. 
Therefore in this case, the map $E_\fin$ is already regarded as a cluster algebraic analogue of the earthquake map. Let us call $E:=E_\fin$ the \emph{cluster earthquake map} in this case. 
Our first result is an analogue of the earthquake theorem for finite type cluster varieties:

\begin{thm}[Earthquake theorem for finite type: \cref{lem:equivariance} and \cref{thm:Eq_finite_type}]\label{introthm:Eq_thm}
For any mutation class $\bs$ of finite type, the cluster earthquake map
\begin{align*}
    E: \cX_\bs(\pos) \times \cX_\bs(\bR^\trop) \to \cX_\bs(\pos)
\end{align*}
is $\Gamma_\bs$-equivariant. Its restriction
\begin{align*}
    E(g_0,-): \X_\bs(\bR^\trop) \xrightarrow{\sim} \X_\bs(\pos) 
\end{align*}
to $\{g_0\} \times \cX_\bs(\bR^\trop)$ is a homeomorphism for any point $g_0$.
\end{thm}
For the rank two cases $A_1 \times A_1,A_2,B_2,G_2$, it is verified by a degree argument 
of the restriction $E(g_0, -)^*: \cX_\bs(\bR^\trop) \setminus \{0\} \to \cX_\bs(\pos) \setminus \{g_0\}$.
A visualization of the image of cluster earthquake map is given in the proof of \cref{lem:homeo_rank2}. 
The general case is proved by an induction on the rank, utilizing the cluster reduction technique \cite{Ish19}. We remark that the finite types $A_n$ and $D_n$ are also surface type, which correspond to $(n+3)$-gon and once-punctured $n$-gon, respectively. 

For a general mutation class $\bs$ of seeds, it seems a highly non-trivial problem to find a correct way to extend the pre-earthquake map to the entire tropical space $\X_\bs(\bR^\trop)$ in a $\Gamma_\bs$-equivariant manner. When $\bs$ is not of finite mutation type, the support $|\fF^+_\bs|$ of the Fock--Goncharov fan is not dense in $\X_\bs(\bR^\trop)$ \cite{Yurikusa}. The case of finite mutation type will be explored in our future work, based on ``cluster twist flows'' generalizing the Fenchel--Nielsen twists. 

\subsection{An application and observations}
Here we give an applications and two observations that shortly follows from \cref{introthm:Eq_thm}. 

\subsubsection{Application: equivariant trivialization of the tangent bundle $T\X_\bs(\pos)$}

The applications of the original earthquake map for closed surfaces is the identification of the tangent bundle $T \cT(\Sigma)$ of the \Teich\ space and the lamination bundle $\ML_\Sigma$ as topological fiber bundles:
\begin{align}\label{eq:tangent_bundle}
    \ML_\Sigma \xrightarrow{\sim} T\cT(\Sigma)
\end{align}
Here $\ML_\Sigma$ is a topological fiber bundle over the \Teich\ space $\cT(\Sigma)$ with the fiber over a marked hyperbolic surface $[X, f] \in \cT(\Sigma)$ being the space $\ML(X)$ of measured geodesic laminations on $X$. 
The isomorphism \eqref{eq:tangent_bundle} is obtained by differentiating the earthquake maps:
\begin{align}\label{eq:quakeflow}
    ([X,f], L) \mapsto \Big([X,f],\ \frac{d}{dt}\bigg|_{t = 0^+} E\big([X,f],\ t \cdot L \big) \Big).
\end{align}
Here $t\cdot$ denotes the $\bR_{>0}$-action on the fibers of $\ML_\Sigma$ that rescales the transverse measure.
The enhanced version for marked surfaces without boundary is also proved (see, for instance, \cite[Proposition 3.36]{BB09}). We give an analogue of this result for cluster varieties of finite type based on the earthquake theorem (\cref{introthm:Eq_thm}): 

\begin{thm}[\cref{thm:tangent_trop}]\label{introthm:tangent_trop}
For any mutation class $\bs$ of finite type, the map
\begin{align*}
    dE: \cX_\bs(\bR_{>0}) \times \cX_\bs(\bR^\trop) \xrightarrow{\sim} T\cX_\bs(\bR_{>0}),\quad
    (g, L) \mapsto \Big(g,\ \frac{d}{dt}\bigg|_{t=0^+} E(g,\, t \cdot L) \Big)
\end{align*}
is a $\Gamma_\bs$-equivariant isomorhpism of topological fiber bundles.
\end{thm}
Notice that since the space $\X_\bs(\pos)$ is contractible, it is obvious that the tangent bundle $T\X_\bs(\pos)$ admits a trivialization. However, the existence of a $\Gamma_\bs$-equivariant trivialization is non-trivial.   

\subsubsection{Observation: deformations of the cluster complex}
As seen from the pictures in the proof of \cref{lem:homeo_rank2}, the cluster pre-earthquake map $E(g_0,-)$ based at $g_0 \in \X_\bs(\pos)$ gives a non-linear deformation of the Fock--Goncharov fan $\fF^+_\bs$. Moreover, looking from a distance, it can be observed that the asymptotic behavior of the curved lines at infinity resembles the structure of the cluster complex of the opposite mutation class $-\bs$ (\cref{def:opposite}). Furthermore, we can vary $g_0$ to get another deformations. We state two results on the asymptotics of $E(g_0,L)$, $(g_0,L) \in \X_\bs(\pos) \times \X_\bs(\bR^\trop)$ along these deformations.

\smallskip
\paragraph{\textbf{Asymptotics as $L$ goes to infinity}}
In order to describe the asymptotic behavior at infinity, let us consider the \emph{Thurston compactification} (a.k.a. the \emph{tropical compactification} or the \emph{Fock--Goncharov compactification}) \cite{FG16,Le16,Ish19}
\begin{align*}
    \overline{\cX_\bs(\pos)} := \cX_\bs(\pos) \sqcup \bS \cX_\bs(\bR^\trop).
\end{align*}
Here $\bS \cX_\bs(\bR^\trop) := (\cX_\bs(\bR^\trop) \setminus \{0\}) / \pos$.
It is endowed with a topology such that a divergent sequence $g_n$ in $\cX_\bs(\pos)$ converges to $[L] \in \bS \cX_\bs(\bR^\trop)$ in $\overline{\cX_\bs(\pos)}$ if and only if $[\log X_{(v)}(g_n)]$ converges to $[x_{(v)}(L)]$ in $\bS \bR^I=(\bR^I \setminus \{0\})/\bR_{>0}$. 
For $L \in \cX_\bs(\bR^\trop)$ and $v \in \bExch_\bs$, we define
\begin{align*}
    \cE^+_{(v)}(g_0) := \Big\{ \lim_{t \to \infty} E(g_0,\, t \cdot L)\ \Big|\ L \in \cC^+_{(v)} \Big\} \subset \bS \cX_\bs(\bR^\trop).
\end{align*}
The following theorem states that the set of the faces of $\{\cE^+_{(v)}(g_0) \mid v \in \bExch_\bs\}$ coincides with the cluster complex \cite{FZ-CA2}
$\bS \fF^+_{-\bs} := \{ \bS \cC \mid \cC \in \fF^+_{-\bs}\}$ for the opposite mutation class $-\bs$:

\begin{thm}[\cref{thm:D_at_bdry}]\label{introthm:asymptotic_L}
Under the canonical identification $\iota: \X_\bs(\bR^\trop) \xrightarrow{\sim} \X_{-\bs}(\bR^\trop)$ (see \eqref{eq:iota}), we have
\begin{align*}
    \cE^+_{(v)}(g_0) = \iota^{-1}(\bS \cC^+_{(-v)}) \subset \bS \cX_\bs(\bR^\trop)
\end{align*}
for each $g_0 \in \cX_\bs(\pos)$ and $v \in \bExch_\bs$.
\end{thm}

\smallskip
\paragraph{\textbf{Asymptotics as $g_0$ goes to infinity}}
Next, we are going to describe the asymptotic behavior of the subsets $\cD_{(v)}^+(g_0)$ in the log-coordinate system $\log X_{(v_0)}(\cX_\bs(\bR^\trop))$ for $v_0 \in \bExch_\bs$ as $g_0$ to the Thurston boundary. 

For $g \in \cX_\bs(\pos)$ and $v \in \bExch_\bs$ we define
\begin{align}\label{eq:shear_MS_cl}
    \mathbf{u}^{(v)}_{g}: \cX_\bs(\bR^\trop) \to \bR^I,\quad
    L \mapsto \log \frac{\mathbf{X}^{(v)}(E(g,L))}{\mathbf{X}^{(v)}(g)}.
\end{align}
This map is nothing but the ``shear coordinates'' in the sense of Meusburger--Scarinci \cite[(17)]{MS}, which is an ingredient of parametrizations of the gravitational moduli space on $\Sigma\times \bR$.

The following theorem states that the subsets 
\begin{align*}
    \mathbf{u}^{(v_0)}_{g}(\cC^+_{(v)}) = \log \frac{\mathbf{X}^{(v_0)}}{\mathbf{X}^{(v_0)}(g_0)}(\cD_{(v)}(g_0))
\end{align*}
converge to $x_{(v_0)}(\cC^+_{(v)})$ in the Gromov--Hausdorff sense as $g_0  \in \cX_\bs(\pos)$ diverges in an appropriate direction:

\begin{thm}[\cref{thm:g_to_infty}]\label{introthm:asymptotic_g}
Let $\bs$ be a mutation class of finite type and $v_0 \in \bExch_\bs$.
Then, 
\begin{align*}
    \mathbf{u}^{(v_0)}_g(\fF^+_\bs) \to \mathbf{x}^{(v_0)}(\fF^+_{\bs})
\end{align*}
as $g_0$ diverges so that $\mathbf{X}^{(v_0)}(g_0)^{-1} = (X^{(v_0)}_i (g_0)^{-1})_{i \in I} \to \mathbf{0}$.
\end{thm}
The proof is based on the positivity result of Gupta on \emph{deformed $F$-polynomials} (\cite[Section 5]{Gup}).

In view of \cref{introthm:asymptotic_L,introthm:asymptotic_g} the simplicial sets $\{\cE^+_{(v)}(g_0)\mid v \in \bExch_\bs\}$ and $\mathbf{u}^{(v_0)}_{g}(\fF_{\bs}^+)$ can be view as ``continuous deformations'' of the cluster complex (or the Fock--Goncharov fan). 
Moreover, although the ``shear coordinates'' $\mathbf{u}^{(v)}_{g}$ are different from our shear coordinates, \cref{introthm:asymptotic_g} tells us that these coordinates are asymptotically the same.

\subsection{Future topics of research}
The above theorems related to the cluster earthquake maps just open up a new direction of research on cluster algebras in the spirit of \Teich--Thurston theory. A clear continuation will be to explore the infinite type cases, which will be discussed elsewhere. Here we present a few related topics. 

\subsubsection{Relation to the space of stability conditions and Mirzakhani's work}
In the case where $\Sigma$ is a closed surface, Mirzakhani \cite{Mir} showed a striking result that the earthquake flow induces an ergodic flow on the space $\ML^1_\Sigma/MC(\Sigma)$, where $\ML^1_\Sigma \subset \ML_\Sigma$ is the unit sub-bundle with respect to the length pairing. Her proof is based on a constructon of a $MC(\Sigma)$-equivariant measurable map $F: \ML^1_\Sigma \to \cQ_\Sigma^1$ that conjugates the earthquake flow to the \emph{\Teich\ horocyclic flow} on the unit sub-bundle $\cQ_\Sigma^1 \subset \cQ_\Sigma$ of the vector bundle of holomorphic quadratic differentials on $\cT(\Sigma)$. Then the ergodicity of earthquake flow follows from that of the \Teich\ horocyclic flow, which has been known before. 

Regarding the relation to the cluster varieties, we remark that the bundle of quadratic differentials is closely related to the Bridgeland's space of \emph{stability conditions} \cite{Bri07} on the finite dimensional derived category of the Ginzburg dg algebra of a quiver with potential \cite{BS}. In a future work, we will discuss an analogue of the \Teich\ horocyclic flow on the space of stability conditions and a measurable conjugacy with our cluster earthquake flow. 

\subsubsection{Higher-rank case: eruption-shearing flows}
Our notion of cluster pre-earthquake naturally includes the \emph{eruption-shearing flows} studied by Sun--Wienhard--Zhang \cite{SWZ20}. The eruption-shearing flows are certain Hamiltonian flows on the Hitchin compoent $\mathrm{Hit}_{n}(\Sigma) \subset \mathrm{Hom}(\pi_1(\Sigma),PSL_n(\bR))/PSL_n(\bR)$ defined as deformations of Frenet curves associated to Hitchin representations, which generalizes the Fenchel--Nielsen twist flows for $n=2$. 
These flows can be also defined on the positive-real part $\X_{PGL_n,\Sigma}(\pos)$ of the moduli space of framed $PGL_n$-local systems \cite{FG03}. The moduli space $\X_{PGL_n,\Sigma}$ admits a natural cluster structure encoded in a mutation class $\bs(\mathfrak{sl}_n,\Sigma)$ such that $\X_{PGL_n,\Sigma}(\pos) \cong \X_{\bs(\mathfrak{sl}_n,\Sigma)}(\pos)$. Then by \cite[Propositions 3.7 and 3.13]{SWZ20}, one can see that the eruption-shearing flows have the form of exponential map (\cref{def:exp}) in the Fock--Goncharov coordinates. Therefore these flows can be regarded as the cluster pre-earthquake flows along certain elementary laminations in $\X_{\bs(\mathfrak{sl}_n,\Sigma)}(\bZ^\trop)$. 

In the case $n=3$, the moduli space $\X_{PGL_3,\Sigma}(\pos)$ parametrizes convex $\bR \bP^2$-structures on $\Sigma$ \cite{FG07c}. The eruption-shearing flows are originally described as explicit deformations of such geometric structures. The tropical space $\X_{\bs(\mathfrak{sl}_3,\Sigma)}(\bZ^\trop)$ is described in terms of Kuperberg's $\mathfrak{sl}_3$-webs in \cite{IK22}. Then the eruption flows that deforms the face coordinates are the cluster pre-earthquake flows along elementary $\mathfrak{sl}_3$-laminations given by a web with unique trivalent vertex and three legs. It will be an interesting problem to define an extension $\X_{\bs(\mathfrak{sl}_3,\Sigma)}(\bR^\trop) \times \X_{\bs(\mathfrak{sl}_3,\Sigma)}(\pos) \to \X_{\bs(\mathfrak{sl}_3,\Sigma)}(\pos)$ in terms of convex $\bR \bP^2$-structures and/or $\mathfrak{sl}_3$-webs. 
We also remark that the mutation class $\bs(\mathfrak{sl}_3,\Sigma)$ is finite type $A_1,D_4,E_7$ when $\Sigma$ is a $3,4,5$-gon respectively. Therefore the cluster nature of the eruption-shearing flows in these cases is contained in our results. 

\subsection*{Organization of the paper}
In \cref{sec:cluster}, basic materials on the cluster algebras are recollected. We introduce the cluster earthquake maps in \cref{sec:cl_eq}, and prove \cref{introthm:Eq_thm}. \cref{introthm:tangent_trop} is proved in \cref{sec:tangent}. Then we proceed to investigate the asymptotic behaviors of earthquake maps in \cref{sec:asymptotics}, and prove \cref{introthm:asymptotic_L,introthm:asymptotic_g}. 

To motivate the constructions, we recall the cluster variety $\X_\Sigma$ associated with a marked surface and their relations to the \Teich/lamination spaces in \cref{cl_surf}. In \cref{sec:eq}, we briefly review Thurston's earthquake map and Bonsante--Krasnov--Schlenker's generalization of it.

\subsection*{Acknowledgements}
The authors thank to Yasuaki Gyoda and Osamu Iyama for informing the existence of green sequences for mutation classes of finite type. Their thanks also goes to Kyoung-Do Kim for valuable discussion on generalized shear coordinates and the work \cite{BB09} of Benedetti--Bonsante. 
T.~A.~would like to thank his supervisor, Takuya Sakasai for his patient guidance.
T.~I.~is partially supported by JSPS KAKENHI (20K22304).
S.~K.~is partially supported by scientific research support of Research Alliance Center for Mathematical Sciences, Tohoku University.

\section{Cluster algebras}\label{sec:cluster}
In this section, we recollect basic notations of \emph{cluster algebras} \cite{FZ-CA4} and \emph{cluster varieties (manifolds)} \cite{FG09}. We partially follow the way of presentation given in \cite{Nak21-1} and \cite{IK20}, respectively. 


\smallskip
Let $I := \{1, 2, \dots, n\}$ and $\mathcal{F}_X\ceq\bQ (X_i^0 \mid i \in I)$ be the field of rational functions on indeterminates 
 $X^0_i\ (i\in I)$.
\subsection{Labeled seeds and their mutations}

\begin{defi}
    A \emph{labeled (skew-symmetrizable) seed} is 
    a pair $(\ve,\mathbf{X})$, where
    \begin{enumerate}
\renewcommand{\labelenumi}{(\roman{enumi})}
        \item $\ve=(\ve_{ij})_{i,j\in I}\in\bZ^{I\times I}$ is a skew-symmetrizable matrix. Namely, there exists positive integers $(d_i)_{i \in I}$ such that $\ve_{ij}d_j = \ve_{ji}d_i$ for all $i,j \in I$;
        \item $\mathbf{X}=(X_i)_{i \in I}$ is a transcendence basis of $\mathcal{F}_X$. Namely, $\cF_X \cong \bQ(X_i \mid i \in I)$.
    \end{enumerate} 
    The matrix $\ve$ is called the \emph{exchange matrix}, and the variables 
    $X_i$ are called 
    the \emph{$\mathcal{X}$-variables}.
\end{defi}
Let $\cS$ denote the set of all labeled seeds in $\cF_X$.
Let $\sgn(x) \in \{-1,0,+1\}$ denote the sign of a real number $x$: $\sgn(x)=-1$ if $x<0$, $\sgn(x)=0$ if $x=0$, and $\sgn(x)=+1$ if $x>0$. 
\begin{defi}
    For $k \in I$, we define the \emph{seed mutation} $\mu_k\cl\cS\to \cS$, $(\ve
        ,\mathbf{X})\mapsto(\ve'
        ,\mathbf{X}')$ at $k \in I$ by 
    \begin{align}
        \ve'_{ij}&=
        \begin{cases}
            -\ve_{i j} & \mbox{if $i=k$ or $j=k$,}\\
            \varepsilon_{i j}+\dfrac{\left|\varepsilon_{i k}\right| \varepsilon_{k j}+\varepsilon_{i k}\left|\varepsilon_{k j}\right|}{2} & \mbox{otherwise},
        \end{cases}\\
        X'_i&=
        \begin{cases}
            X_{k}^{-1} & \mbox{if $i=k$},\\
            X_{i}\big(1+X_{k}^{-\sgn (\ve_{i k})}\big)^{-\ve_{i k}} & \mbox{if $i\neq k$}.
        \end{cases}\label{eq:cluster_transf_X}
    \end{align}
Here $\ve=(\ve_{ij})_{i,j \in I},\ \ve'=(\ve'_{ij})_{i,j \in I},\ 
    \mathbf{X}=(X_i)_{i \in I}\text{ and }  \mathbf{X}'=(X'_i)_{i \in I}$.
The transformation $\mathbf{X} \mapsto \mathbf{X}'$ is called the \emph{cluster $\cX$-transformation} at $k$.
\end{defi}
It is easy to verify that the seed mutation $\mu_k$ is involutive. 
Let $\fS_I^{\mathrm{cl}} \subset \fS_I$ be the subgroup consisting of permuations $\sigma$ such that $d_{\sigma^{-1}(i)}=d_i$ for all $i \in I$. For $\sigma \in \fS_I^{\mathrm{cl}}$, we similarly define $\sigma\cl\cS \to \cS$, $(\ve
    ,\mathbf{X})\mapsto(\ve'
    ,\mathbf{X}')$
by 
\begin{align*}
    \ve'_{ij}=\ve_{\sigma^{-1}(i),\sigma^{-1}(j)} \quad \mbox{and} \quad 
    X'_i= X_{\sigma^{-1}(i)}. 
\end{align*}
Then $\ve'$ is again skew-symmetrizable with the same symmetrizer $(d_i)_{i \in I}$.
\begin{defi}
    We say that two labeled seeds $ (\ve
    , \mathbf{X}),(\ve'
    ,\mathbf{X}') \in \cS$ are \emph{mutation-equivalent} if there is a finite composition of seed mutations and permutations that maps $ (\ve
    , \mathbf{X})$ to $(\ve'
    ,\mathbf{X}')$. A mutation-equivalence class $\bs$ of labeled seeds is simply called a \emph{mutation class}. 
\end{defi}

Mutation classes of labeled seeds are the basic subjects in the research field of cluster algebra. 

\begin{defi}
The relations among the labeled seeds in a given mutation class $\bs$ can be encoded in the \emph{(labeled) exchange graph} $\bExch_\bs
$. It is a graph with vertices $v$ corresponding to the labeled seeds $\bs^{(v)}$ in $\bs$, together with labeled edges of the following two types:

\begin{itemize}
    \item labeled edges of the form $v \overbar{k} v'$ whenever the seeds $\bs^{(v)}$ and $\bs^{(v')}$ are related by the mutation $\mu_k$ for some $k \in I$;
    \item labeled edges of the form $v \overbarnear{\sigma} v'$ whenever the seeds $\bs^{(v)}$ and $\bs^{(v')}$ are related by a transposition $\sigma=(j\ k) \in \fS_I^{\mathrm{cl}}$.
\end{itemize}

\end{defi}

When no confusion can occur, we simply denote a vertex of the labeled exchange graph by $v \in \bExch_\bs$ instead of $v \in V(\bExch_\bs)$. For each vertex $v \in \bExch_\bs$, we denote the corresponding labeled seed by $\bs^{(v)}=(\ve^{(v)},\mathbf{X}^{(v)})$, $\ve^{(v)}=(\ve_{ij}^{(v)})_{i,j \in I}$ and $\mathbf{X}^{(v)}=(X_i^{(v)})_{i \in I}$. 

    Given a skew-symmetrizable matrix $\ve\in\bZ^{I\times I}$, there exists 
    a unique mutation class $\bs$ 
    containing the labeled seed $(\ve,
    (X_i^0)_{i \in I})$.
    We call $(\ve, (X_i^0)_{i \in I})$ the \emph{initial seed} in this mutation class.

\subsection{Classification of mutation classes}

A classification of mutation classes is spelled out in the following sense.

\begin{defi}\label{cl_isom}
    Two mutation classes $\bs_1$ and $\bs_2$ of labeled seeds in $\cF_X$ are \emph{isomorphic} if there exist a field automorphism $f\cl \cF_X\rightarrow\cF_X$ and a graph isomorphism $\alpha\cl \bExch_{\bs_1}\rightarrow\bExch_{\bs_2}$ such that $f(X_i^{(v)})=X_i^{\alpha(v)}$ for any $i\in I$ and any $v\in\bExch_{\bs_1}$. 
     
\end{defi}
\begin{rem}
    The isormorphism defined above is the $\X$-variable version of the \emph{strong isomorphism} of cluster algebras \cite{FZ-CA2}. 
    These notions are indeed the same, thanks to the synchronicity theorem \cite{Nak21-4} (see also \cite[Theorem 2.56]{Ish20}).
    \end{rem}

\begin{defi}A mutation class $\bs$ is said to be
    \begin{enumerate}
        \item of \emph{finite type} if the number of elements in $\bs$ is finite.
        \item \emph{decomposable} if any exchange matrix $\varepsilon$ in $\bs$ is decomposable into a direct sum of two square matrices. Otherwise \emph{indecomposable}. 
    \end{enumerate}
\end{defi}
Let $C=(C_{ij})_{i,j \in I}$ be the Cartan matrix of Dynkin type $Z_n$, $Z \in \{A,B,C,D,E,F,G\}$. Let $\ve$ be any skew-symmetrizable matrix satisfying 
\begin{align}\label{eq:Dynkin_quiver}
    |\ve_{ij}|=-C_{ij}
\end{align}
for $i \neq j$. Then the mutation class containing the seed $(\ve,(X_i^0)_{i \in I})$ is denoted by $Z_n$.

\begin{thm}{\cite[Theorem 1.4]{FZ-CA2}}\label{thm:finite_seed_pattern}
    Any indecomposable mutation class of finite type is either isomorphic to $A_n(n\geq 1),\, B_n(n\geq 2),\, C_n(n\geq 3),\, D_n(n\geq 4),\, E_6, E_7, E_8, F_4, G_2$.
\end{thm}
In particular, mutation classes of finite type of rank $2$ are $A_1 \times A_1$, $A_2$, $B_2$ and $G_2$.

\subsection{Cluster manifolds and their tropical analogues}
Let $\bs$ be a mutation class of labeled seeds in $\cF_X$. 

\begin{defi}\label{def:cluster_manifold}
The \emph{cluster $\cX$-manifold} $\cX_\bs(\pos)$ associated with $\bs$ is a $C^\omega$-manifold homeomorphic to $\pos^I$, equipped with an atlas consisting of global charts $\mathbf{X}^{(v)}: \cX_\bs(\pos) \xrightarrow{\sim} \pos^I$ parametrized by the vertices $v \in \bExch_\bs$ such that the coordinate transformations among them are given by the corresponding cluster $\cX$-transformations and permutations.
\end{defi}

We also have the following tropical analogue of the cluster $\X$-manifold.
Let $f(X_1, \dots, X_n)$ be a positive rational function on $n$ variables (namely, a rational function admitting a subtraction-free expression). 
Then its tropical limit $f^\trop(x_1, \dots, x_n)$ is defined by
\begin{align}\label{eq:trop_limit}
    f^\trop(x_1, \dots, x_n) 
    := \lim_{\epsilon \to -0} \epsilon \log f(e^{x_1/\epsilon}, \dots, e^{x_n/\epsilon}).
\end{align}
One can verify that it defines a piecewise linear function on $\bR^I$. Observe that the cluster $\X$-transformation formula \eqref{eq:cluster_transf_X} is an $n$-component positive rational function. Its tropical limit $\mathbf{x} \mapsto \mathbf{x}'$, which is explicitly given by
\begin{align*}
    x'_i =
    \begin{cases}
        -x_k & \mbox{if } i = k,\\
        x_i - \ve_{ik} \min\{0, -\sgn(\ve_{ik}) x_k\} & \mbox{if } i \neq k,
    \end{cases}
\end{align*}
is called the \emph{tropical cluster $\X$-transformation}.

\begin{defi}\label{def:X^trop}
The \emph{tropical $\X$-manifold} $\X_\bs(\bR^\trop)$ associated with $\bs$ is a piecewise linear manifold homeomorphic to $\bR^I$, equipped with an atlas consisting of global charts $\mathbf{x}^{(v)}: \cX_\bs(\bR^\trop) \xrightarrow{\sim} \bR^I$ parametrized by the vertices $v \in \bExch_\bs$ such that the coordinate transformations among them are given by the tropical cluster $\cX$-transformations and permutations.
\end{defi}

We note that there is an $\bR_{>0}$-action on $\cX_\bs(\bR^\trop)$ so that $\mathbf{x}^{(v)}(t \cdot L) = t \mathbf{x}^{(v)}(L)$ for any $L \in \cX_\bs(\bR^\trop)$, $v \in \bExch_\bs$ and $t \in \bR_{>0}$.

\begin{rmk}
In \cref{def:X^trop}, $\bR^\trop$ denotes the set of real numbers with a structure of semifield such that the addition given by minimum and the multiplication given by the usual addition.
(It is known as \emph{min-plus semifield}.)
The tropical limit \eqref{eq:trop_limit} is equivalently obtained by applying the formal replacements $+ \mapsto \min$ and $\times \mapsto +$ for a subtraction-free expression of the rational function.
More precise and general definition of the tropicalization is given as the $\bR^\trop$-valued points of the positive scheme.
We refer the reader to \cite[Section 2]{GHKK} for more details.
\end{rmk}

\begin{defi}
The \emph{cluster modular group} associated with $\bs$ is the subgroup $\Gamma_\bs \subset \mathrm{Aut}(\bExch_\bs)$ consisting of the graph automorphisms which preserves the exchange matrices assigned to the vertices and the labels on the edges. The cluster modular group acts on $\X_\bs(\pos)$ as $C^\omega$-diffeomorphisms so that $X_i^{(v)}(\phi(g))=X_i^{(\phi^{-1}(v))}(g)$ for all $\phi\in \Gamma_\bs$, $g \in \X_\bs(\pos)$ and $v \in \bExch_\bs$. It acts on $\X_\bs(\bR^\trop)$ as PL automorphisms in the same manner. 
\end{defi}

\subsection{$C,G,F$-matrices and the separation formula}

We use the notation $[a]_+:=\max\{a,0\}$ for a real number $a \in \bR$. 
Let $\bs$ be a mutation class of labeled seeds in $\cF_X$. 
\begin{defi}\label{c_mat}
    Fix a vertex $v_0 \in \bExch_\bs$.
    
    \begin{enumerate}
        \item We assign to each vertex $v$ of $\bExch_\bs$ an integral matrix $C^\bs_{v_0 \to v}=(c_{ij}^{(v)})_{i,j\in I}$
        with the initial condition $C^\bs_{v_0 \to v_0} := \mathrm{Id}$ and the relations
        \begin{align*}
            &c_{i j}^{(v')} =
            \begin{cases}
            -c_{k j}^{(v)} & \mbox{if } i=k, \\
            c_{i j}^{(v)} + \big[\ve_{jk}^{(v)} \big]_{+} c_{ki}^{(v)} + \ve_{jk}^{(v)} \big[-c_{ki}^{(v)}\big]_{+} & \mbox{if } i \neq k, \end{cases}
        \end{align*}
        for a horizontal edge $v \overbar{k} v'$, and 
        \begin{align*}&c_{i j}^{(v')} =
            c^{(v)}_{\sigma^{-1}(i), j}  
        \end{align*}
        for a vertical edge $v \overbarnear{\sigma} v'$. 
        The matrices $C^\bs_{v_0 \to v}$ are called the \emph{$C$-matrices} from $v_0$ to $v$ in $\bs$, and their row vectors are called the \emph{$c$-vectors} from $v_0$ to $v$ in $\bs$.
        \item \label{item:F_mutation} We assign to each vertex $v$ of $\bExch_\bs$ a tuple $(F_i^{(v)})_{i\in I}$ of polynomials in variables $y_1, \dots, y_n$
        with the initial condition $F^{v_0 \to v_0}_i := 1$ for all $i \in I$ and the relations $F_{i}^{v_0 \to v'} = F_{i}^{v_0 \to v}$, $i \neq k$ and
        \begin{align*}
        \qquad F_{k}^{v_0 \to v'} \cdot F_{k}^{v_0 \to v}
        = \dsp\prod_{j \in I} y_{j}^{[c_{k j}^{(v)}]_+} \prod_{\ell \in I}\big(F_{\ell}^{{v_0 \to v}}\big)^{[\ve_{k \ell}^{(v)}]_+}+\prod_{j \in I} y_{j}^{[-c_{k j}^{(v)}]_+} \prod_{\ell \in I}\big(F_{\ell}^{{v_0 \to v}}\big)^{[-\ve_{k \ell}^{(v)}]_+}
        \end{align*}
        for a horizontal edge $v \overbar{k} v'$, and
        \begin{align*}
            F_{i}^{v_0 \to v'} = F_{\sigma^{-1}(i)}^{v_0 \to v}
        \end{align*}
    for a vertical edge $v \overbarnear{\sigma} v'$.
    We call $F^{v_i \to v}_i$ the $i$-th \emph{$F$-polynomial}\footnote{It is non-trivial from our definition that the $F$-polynomials are actually polynomials. This fact is a consequence of the Laurent phenomenon property of the cluster $\cA$-variables. See \cite{FZ-CA4} for more details.} from $v_0$ to $v$ in $\bs$.
    \end{enumerate}
\end{defi}

\begin{thm}[separation formula {\cite[Proposition 3.13 and Corollary 6.3]{FZ-CA4}}]\label{thm:sep_for}
    For a mutation class $\bs$ with an initial seed $(\ve^0, (X^0_i)_{i \in I})$ at $v_0 \in \bExch_\bs$, 
    we have
    \[
        X_{i}^{(v)}=\prod_{j\in I} (X^0_{j})^{c_{i j}^{(v)}} \cdot  F_{j}^{v_0 \to v}(X_1^0, \dots, X_n^0)^{\ve_{ij}^{(v)}},
    \]
    where $(c_{ij}^{(v)}) = C^\bs_{v_0 \to v}$.
\end{thm}

\begin{defi}\label{def:G-mat}
Fix a vertex $v_0 \in \bExch_\bs$. 
For each $v\in \bExch_\bs$, we define
the \emph{$G$-matrix $G^\bs_{v_0 \to v}$ from $v_0$ to $v$} in $\bs$ is defined by
\begin{align*}
    G^\bs_{v_0 \to v} := D ((C^\bs_{v_0 \to v})^{-1})^\tr D^{-1}.
\end{align*}
Here, $D=\mathrm{diag}(d_i \mid i \in I)$ consists of the symmetrizer of the exchange matrices $\ve^{(v)}$.
\end{defi}

\begin{defi}\label{def:opposite}
For a mutation class $\bs$ with an initial seed $(\ve^0, (X_i)_{i \in I})$, its \emph{opposite mutation class} $-\bs$ is defined to be the 
mutation class containing $(-\ve^0, (X_i^{-1})_{i \in I})$.
\end{defi}
Since $\mu_k(-\ve, \mathbf{X}^{-1})=(-\ve',(\mathbf{X}')^{-1})$ for $(\ve',\mathbf{X}')=\mu_k(\ve,\mathbf{X})$, the mutation class $-\bs$ actually does not depend on the choice of an initial seed.  
The labeled exchange graphs of $\bs$ and $-\bs$ are isomorphic. 
Slightly abusing the notation, 
we denote the labeled seeds in $-\bs$ by 
$\bs^{(-v)} = (\ve^{(-v)}, \mathbf{X}^{(-v)})$ so that $\ve^{(-v)} = - \ve^{(v)}$ and $\mathbf{X}^{(-v)} = (\mathbf{X}^{(v)})^{-1}$ for $v\in \bExch_\bs$. Namely, the corresponding vertex in $\bExch_{-\bs}$ is denoted by $-v$. 

\begin{thm}[{\cite{GHKK} and \cite[Theorem 1.2]{NZ}}]\label{thm:C-mat}
Let $\bs$ be a mutation class, and fix $v_0 \in \bExch_\bs$.
\begin{enumerate}
    \item \label{item:sign_coh} For each $v \in \bExch_\bs$ and $k \in I$, there is a sign $\epsilon \in \{+, -\}$ such that $\epsilon \mathbf{c}_k \in \bZ_{\geq 0}^I$.
    Here, $\mathbf{c}_k$ is the $k$-th $c$-vector from $v_0$ to $v$ in $\bs$.
    \item \label{item:trop_dual} For any $v \in \bExch_\bs$, $C^\bs_{v \to v_0} = (C^{-\bs}_{v_0 \to v})^{-1}$.
\end{enumerate}
\end{thm}
The phenomenon \eqref{item:sign_coh} is called the \emph{sign coherence} of $c$-vectors.
We will refer to the sign $\epsilon$ as the \emph{$k$-th tropical sign at $v$}.
Our definition of $G$-matrices in \cref{def:G-mat} is equivalent to the original definition given
in \cite{FZ-CA4}, thanks to \cite{NZ}. 
The relations given in \cref{def:G-mat} and \cref{thm:C-mat} \eqref{item:trop_dual} are referred to as the \emph{tropical duality} of $C$- and $G$-matrices.

\subsection{Fock--Goncharov fan}
Let $\bs$ be a mutation class of labeled seeds in $\cF_X$.
\begin{defi}\label{def:positive_cone}
For $v \in \bExch_\bs$,
we define
\begin{align*}
    \cC_{(v)}^{+}:= \{ w \in \X_\bs(\bR^\trop) \mid
    x_i^{(v)}(w) \geq 0,\ \text{for any }i \in I \},
\end{align*}
\emph{i.e.}, the cone of points with non-negative coordinates at $v$.
Let $\fF_\bs^+$ denote the set of the faces of the cones $\cC^+_{(v)}$, $v \in \bExch_\bs$. We call it the \emph{Fock--Goncharov fan} associated with $\bs$.
We denote by $|\fF^+_\bs|$ the union $\bigcup_{v \in \bExch_\bs}\cC^+_{(v)}$.
\end{defi}

\begin{thm}{\cite[Corollary 5.5 and Theorem 5.8]{FG09}}\label{complete_fan}
Let $\bs$ be any mutation class.
For each $v_0 \in \bExch_\bs$,
\begin{enumerate}
    \item The set $\mathbf{x}^{(v_0)}(\fF^+_\bs) = \{ \mathbf{x}^{(v_0)}(\cC) \mid \cC \in \fF^+_\bs \}$ is a fan in $\bR^I$.
    \item Moreover if $\bs$ is of finite type, then the fan $\mathbf{x}^{(v_0)}(\fF^+_\bs)$ is complete.
    Namely, $\mathbf{x}^{(v_0)}(|\fF^+_\bs|) = \bR^I$.
\end{enumerate}
\end{thm}

As the following proposition states, the cones in $\fF_\bs^+$ is spanned by $g$-vectors. So they are also called the \emph{$g$-vector fan}.

\begin{prop}\label{prop:cl_cpx=g-vect_fan}
Fix $v_0 \in \bExch_\bs$.
For $v \in \bExch_\bs$, the cone $\mathbf{x}^{(v_0)}(\cC^+_{(v)})$ is spanned by the $g$-vectors from $v_0$ to $v$ in $-\bs^\vee$.
Here, $\bs^\vee$ is the mutation class with an initial seed $(D \ve^{(v_0)} D^{-1}, ((X^{(v_0)}_i)^{d_i})_{i \in I})$. 
\end{prop}
\begin{proof}
By \cite{IK21}, the presentation matrix of the linear extension of
\begin{align*}
    \mathbf{x}^{(v_0)} \circ (\mathbf{x}^{(v)})^{-1}:\mathbf{x}^{(v)}(\cC^+_{(v)}) = \bR_{\geq 0}^I \to \mathbf{x}^{(v_0)}(\cC^+_{(v)}) \subset \bR^I
\end{align*}
is the $C$-matrix $C^{\bs}_{v \to v_0}$.
By \cref{thm:C-mat} \eqref{item:trop_dual}, we have
\begin{align*}
    \mathbf{x}^{(v_0)} \circ (\mathbf{x}^{(v)})^{-1}(\mathbf{e}_i)
    &= (\mbox{$i$-th column vector of $C^{\bs}_{v \to v_0}$})\\
    &= (\mbox{$i$-th column vector of $(C^{\bs}_{v_0 \to v})^{-1}$})\\
    &= (\mbox{$i$-th row vector of $D G^{-\bs}_{v_0 \to v} D^{-1}$}).
\end{align*}
From the definition, a row vector of $G^{-\bs}_{v_0 \to v}$ is the $g$-vector in $-\bs^\vee$.
\end{proof}

\section{Cluster earthquake maps}\label{sec:cl_eq}
Our ultimate goal is to define a cluster algebraic analogue
\begin{align*}
    E\cl \cX_\bs(\bR_{>0})\times\cX_\bs(\bR^\trop) \to \cX_\bs(\bR_{>0})
\end{align*}
of earthquake maps for any mutation classes not coming from marked surfaces.
In the surface case, $\X_\bs(\pos)$ is the enhanced \Teich\ space and $\X_\bs(\bR^\trop)$ is a space of measured geodesic laminations (see \cref{cl_surf}). 
As we see in \cref{subsec:eq_exp}, the coordinate expression of the earthquake $E(g_0,L)$ 
takes the form of a multiplication of the exponentials of the shear coodinates of the lamination $L$, as long as it belongs to the Fock--Goncharov fan. 
In view of this fact, we introduce a cluster algebraic analogue of earthquakes and prove the cluster earthquake theorem for mutation classes of finite type.

\subsection{Exponential maps}
Let $\bs$ be any mutation class. 

Fixing a basepoint $g_0 \in \X_\bs(\bR_{>0})$, let us consider the subset
\begin{align}\label{def_D}
    \cD_{(v)}^+(g_0):= \{ g \in \X_{\bs}(\bR_{>0}) \mid X_i^{(v)}(g) \geq X_i^{(v)}(g_0)\ \text{for any }i \in I\}
\end{align}
for each $v \in \bExch_\bs$. On the other hand, recall the cone $\cC^+_{(v)}$ in the tropical $\X$-manifold from \cref{def:positive_cone}. 
Then we define a map 
\begin{align*}
    \exp_{g_0}^{(v)} : \cC_{(v)}^+ \to \cD_{(v)}^+(g_0) \subset \X_\bs(\bR_{>0})
\end{align*}
by 
\begin{align}\label{eq:def_exponential}
    X_i^{(v)} (\exp_{g_0}^{(v)}(L)) := \exp(x_i^{(v)}(L)) \cdot X_i^{(v)}(g_0)
\end{align}
for $L \in \cC^+_{(v)}$ and $i \in I$. 

\begin{lem}
For a horizontal edge $v_1 \overbar{k} v_2$ in $\bExch_\bs$,
we have
\begin{align*}
    \exp_{g_0}^{(v_1)}(\cC^+_{(v_1)} \cap \cC^+_{(v_2)}) = 
    \exp_{g_0}^{(v_2)}(\cC^+_{(v_1)} \cap \cC^+_{(v_2)}) 
\end{align*}
\end{lem}

\begin{proof}
Suppose that $v_1$ and $v_2$ are related by the mutation directed to $k \in I$. Then for each $L \in \cC^+_{(v_1)}\cap \cC^+_{(v_2)}$, we have $x_k^{(v_1)}(L)=-x_k^{(v_2)}(L)=0$ and hence $x_i^{(v_1)}(L)=x_i^{(v_2)}(L)$ for $i \neq k$ from the tropical cluster transformation formula. It follows that
\begin{align*}
    X_k^{(v_1)}(\exp_{g_0}^{(v_2)}(L)) &= X_k^{(v_2)}(\exp_{g_0}^{(v_2)}(L))^{-1} \\
    &= (X_k^{(v_2)}(g_0))^{-1} \\
    &= X_k^{(v_1)}(g_0),
\end{align*}
and 
\begin{align*}
    X_i^{(v_1)}(\exp_{g_0}^{(v_2)}(L)) &= X_i^{(v_2)}(\exp_{g_0}^{(v_2)}(L))\Big( 1+ (X_k^{(v_2)}(\exp_{g_0}^{(v_2)}(L)))^{-\sgn \ve_{ik}^{(v_2)}}\Big)^{-\ve_{ik}^{(v_2)}} \\
    &= e^{x_i^{(v_2)}(L)}X_i^{(v_2)}(g_0)\Big( 1+ (X_k^{(v_2)}(g_0))^{-\sgn \ve_{ik}^{(v_2)}}\Big)^{-\ve_{ik}^{(v_2)}} \\
    &= e^{x_i^{(v_1)}(L)} X_i^{(v_1)}(g_0)\Big( 1+ (X_k^{(v_1)}(g_0))^{-\sgn \ve_{ik}^{(v_1)}}\Big)^{-\ve_{ik}^{(v_1)}}
    \Big( 1+ (X_k^{(v_1)}(g_0))^{-\sgn \ve_{ik}^{(v_1)}}\Big)^{\ve_{ik}^{(v_1)}}\\
    &= e^{x_i^{(v_1)}(L)}X_i^{(v_1)}(g_0)
\end{align*}
for $i \neq k$. Hence the point $\exp_{g_0}^{(v_2)}(L)$ satisfies the same coordinate charaterization as $\exp_{g_0}^{(v_1)}(L)$, and thus $\exp_{g_0}^{(v_2)}(L)=\exp_{g_0}^{(v_1)}(L)$.
\end{proof}

\begin{defi}\label{def:exp}
We define the
\emph{exponential map} $\exp_{g_0}: |\fF^+_\bs| \to \cX_\bs(\bR_{>0})$ by gluing the maps $\exp_{g_0}^{(v)}$ :
\begin{equation*}
    \begin{tikzcd}[row sep=small]
    {|\fF^+_\bs|} \ar[r, "\exp_{g_0}"] \ar[d, phantom, "\cup"] & \cX_\bs(\bR_{>0}) \ar[d, phantom, "\cup"]\\
    \cC_{(v)}^+ \ar[r, "\exp_{g_0}^{(v)}"] & \cD_{(v)}^+(g_0)
    \end{tikzcd}
\end{equation*}
We also define a \emph{cluster pre-earthquake map} \[
    E_\mathrm{pre}\cl  \X_\bs(\pos) \times |\fF_\bs^+| \to \X_\bs(\pos),\ (g_0,L)\mapsto \exp_{g_0}(L).
\]
\end{defi}

\begin{rem}
    For a mutation class $\bs$ of surface type, the cluster pre-earthquake map coincides with the restriction of the earthquake map to $ \X_\bs(\pos) \times |\fF^+_\bs|$ by \cref{eq_surf}. Therefore, a cluster pre-earthquake map is a partial generalization of the earthquake map.
\end{rem}

\begin{lem}\label{lem:equivariance}
A cluster pre-earthquake map $E_\fin$
is $\Gamma_\bs$-equivariant. Namely, for any $(g_0,L)\in\X_\bs(\pos) \times |\fF^+_\bs|$ and  $\phi\in\Gamma_\bs$, we have
\[
    \phi(E_\fin(g_0,L))=E_\fin(\phi(g_0),\phi(L))
\]
\end{lem}

\begin{proof}
Recall that a mutation loop $\phi \in \Gamma_\bs$ acts on $\X_\bs(\pos)$ so that $X_i^{(v)}(\phi(g))=X_i^{(\phi^{-1}(v))}(g)$ for all $g \in \X_\bs(\pos)$ and $v \in \bExch_\bs$. The PL action on $\X_\bs(\bR^\trop)$ is similar. We have $\phi(\cC_{(v)}^+) = \cC_{(\phi(v))}^+$ for all $v \in \bExch_\bs$, since $x_i^{(\phi(v))}(\phi(L))=x_i^{(v)}(L) \geq 0$ for $L \in \cC_{(v)}^+$. 

Fix $v \in \bExch_\bs$, and let $v':=\phi(v)$. 
Then for $L \in \cC_{(v)}^+$, we have $\phi(L) \in \cC_{(v')}^+$ and hence
\begin{align*}
    X_i^{(v')}(\exp_{\phi(g_0)}(\phi(L))) &= \exp(x_i^{(v')}(\phi(L)))\cdot X_i^{(v')}(\phi(g_0)) \\
    &=\exp(x_i^{(v)}(L))\cdot X_i^{(v)}(g_0) \\
    &=X_i^{(v)}(\exp_{g_0}(L)) \\
    &=X_i^{(v')}(\phi(\exp_{g_0}(L))) 
\end{align*}
for all $i \in I$. Thus we get $\exp_{\phi(g_0)}(\phi(L)) = \phi(\exp_{g_0}(L))$ as they have the same coordinates.
\end{proof}

In general, it is highly non-trivial how to extend $\exp_{g_0}$ to a map from the entire space $\X_\bs(\bR^\trop)$: this point will be pursued in future works. In this paper, we discuss mutation classes of finite type, where $\exp_{g_0}$ is sufficient for our purpose.

\subsection{Finite type}
A specific feature of finite type mutation classes is that the Fock--Goncharov fan is complete by Theorem \ref{complete_fan} (2). Therefore, we may simply define the \emph{cluster earthquake map} to be 
\begin{align*}
    E:=E_\fin: \cX_\bs(\pos) \times \X_\bs(\bR^\trop) \to \X_\bs(\pos),\ (g_0,L) \mapsto \exp_{g_0}(L)
\end{align*}
in this case. The following gives a cluster algebraic analogue of the Thurston's earthquake theorem:

\begin{thm}\label{thm:Eq_finite_type}
For any mutation class $\bs$ of finite type,
the cluster earthquake map gives a 
homeomorphism
\begin{align*}
    E(g_0,-): \X_\bs(\bR^\trop) \xrightarrow{\sim} \X_\bs(\bR_{>0})
\end{align*}
for any basepoint $g_0 \in \X_\bs(\bR_{>0})$.
\end{thm}

Below we give a proof of this theorem.

\begin{lem}\label{lem:contin_closed}
For any mutation class $\bs$ of finite type and any basepoint $g_0 \in \X_\bs(\bR_{>0})$, the map $E(g_0,-)$ is a continuous, proper closed map.
\end{lem}

\begin{proof}
It follows from the fact that the restriction $E(g_0,-)|_{\cC_{(v)}^+}=\exp_{g_0}^{(v)}$ to each maximal dimensional cone $\cC_{(v)}^+$ for $v \in \bExch_\bs$ is a homeomorphism onto its image, and that the tropical space $\X_\bs(\bR^\trop)$ is covered by a finite number of such closed sets.
\end{proof}

We are going to proceed by induction on the rank $N \geq 1$ of the mutation class $\bs$. The rank one case is obvious.

\begin{lem}\label{lem:homeo_rank2}
For any rank two mutation classes of finite type, \cref{thm:Eq_finite_type} holds true.
\end{lem}

\begin{figure}
\centering
\begin{tikzcd}[row sep=small]
&\mathbf{x}^{(v_0)}(\cX_{\bs}(\bR^\trop))
\ar[rrr, "{\log \mathbf{X}^{(v_0)} \circ E(g_0,-) \circ (\mathbf{x}^{(v_0)})^{-1}}"] &&& \log \mathbf{X}^{(v_0)}(\cX_{\bs}(\pos))\\
[-3mm]\bs = A_1 \times A_1: \hspace{-5mm} &
{\includegraphics[width=4cm, align=c]{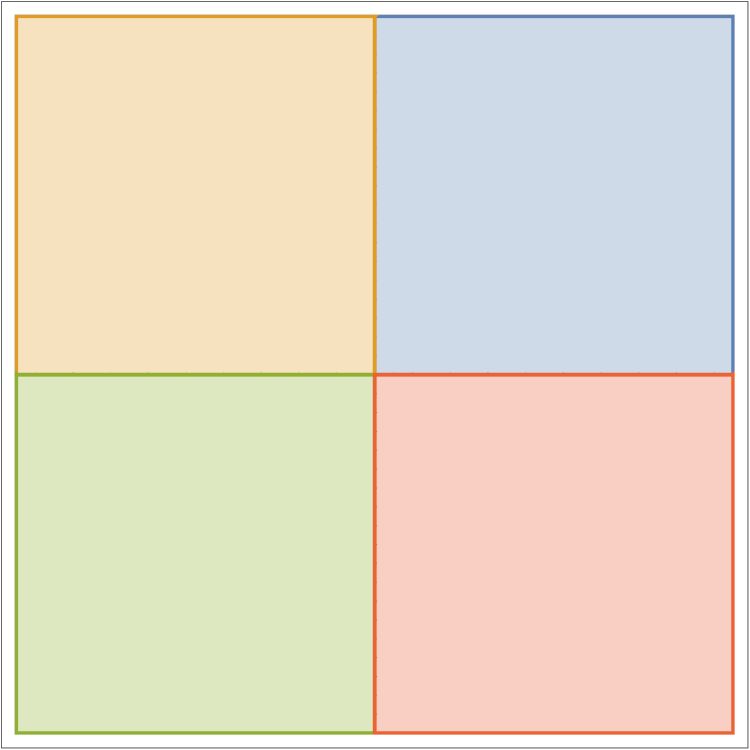}} &&&{\includegraphics[width=4cm, align=c]{A1A1_C.pdf}}\\
\bs = A_2: \hspace{-5mm} & {\includegraphics[width=4cm, align=c]{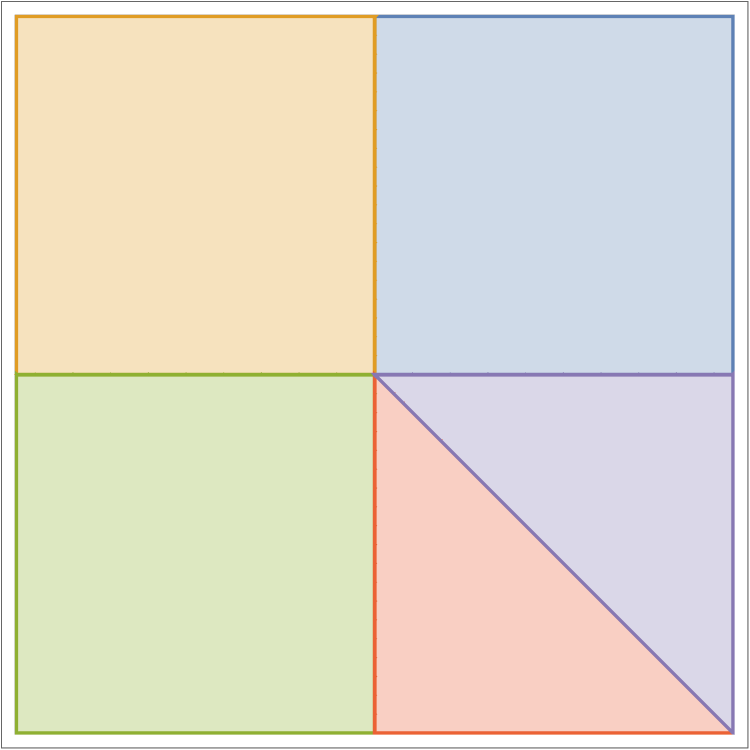}} &&& {\includegraphics[width=4cm, align=c]{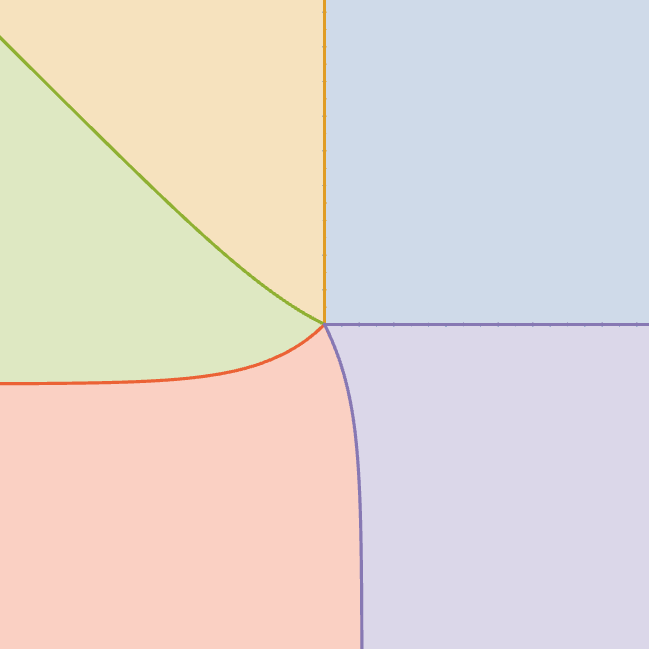}}\\
\bs = B_2: \hspace{-5mm} & {\includegraphics[width=4cm, align=c]{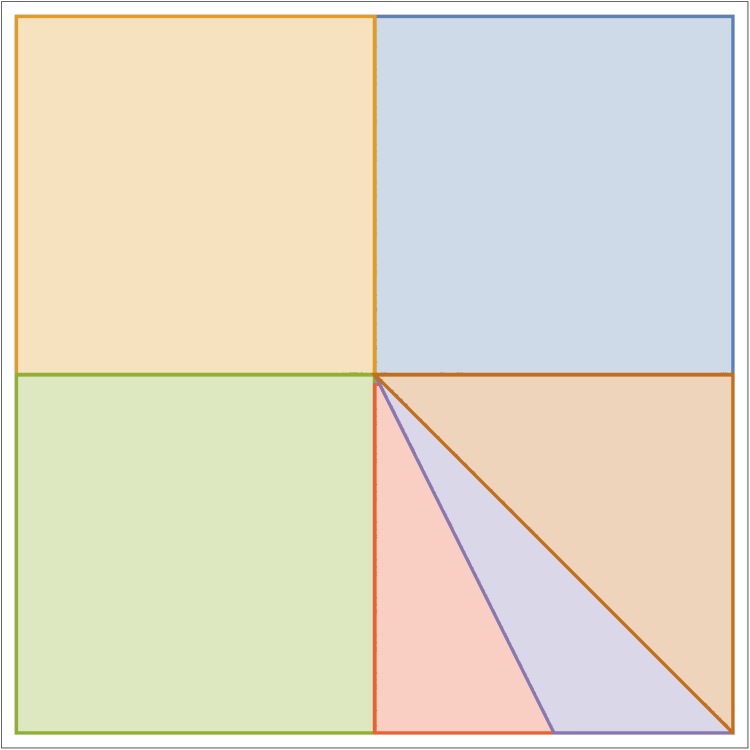}} &&& {\includegraphics[width=4cm, align=c]{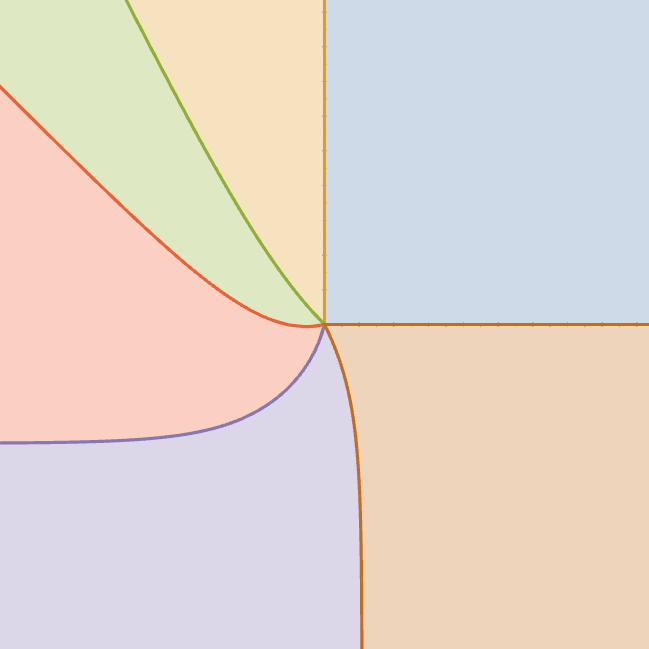}}\\
\bs = G_2: \hspace{-5mm} & {\includegraphics[width=4cm, align=c]{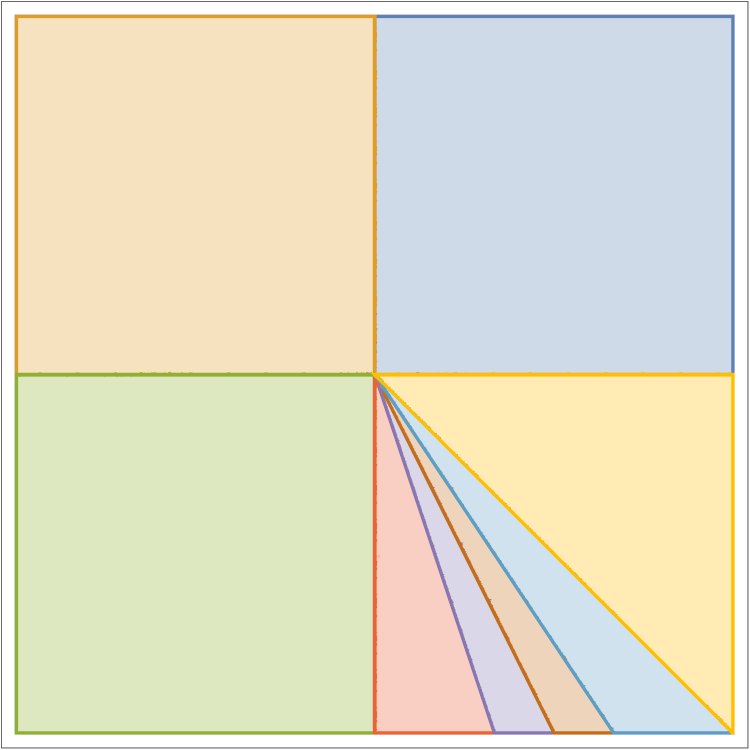}} &&& {\includegraphics[width=4cm, align=c]{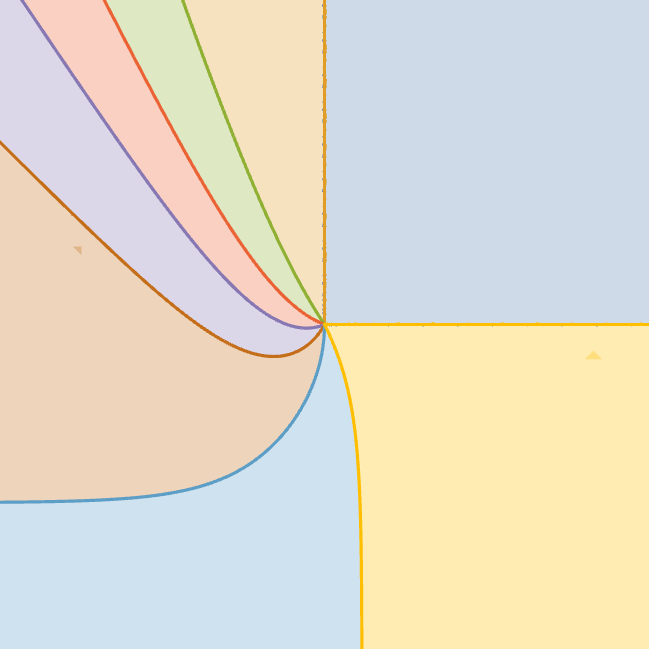}}
\end{tikzcd}
\caption{The pictures of the cones $\mathbf{x}^{(v_0)}(\cC^+_{(v)})$ and the subsets $\log \mathbf{X}^{(v_0)}(\cD^+_{(v)}(g_0))$ of the seed $\bs$ of rank 2.
The same color means the same vertex $v \in \bExch_\bs$.
The exchange matrices at $v_0$ are $\begin{psmallmatrix}0 & 0\\ 0 & 0\end{psmallmatrix}$, $\begin{psmallmatrix}0 & -1\\ 1 & 0\end{psmallmatrix}$, $\begin{psmallmatrix}0 & -1\\ 2 & 0\end{psmallmatrix}$ and $\begin{psmallmatrix}0 & -1\\ 3 & 0\end{psmallmatrix}$ for $\bs = A_1 \times A_2$, $A_2$, $B_2$ and $G_2$, respectively.
These pictures are drown in the range $[-6, 6] \times [-6, 6]$.}
\label{cones}
\label{fig:rk2_eq_plots}
\end{figure}

\begin{proof}
By \cref{thm:finite_seed_pattern}, it suffices to consider the cases of $A_1 \times A_1$, $A_2$, $B_2$ and $G_2$. When we fix a point $g_0\in\cX_\bs(\bR_{>0})$, the restriction $E(g_0,-)^*\cl\cX_\bs(\bR^\trop) \setminus \{0\}\rightarrow\cX_\bs(\bR_{>0}) \setminus \{g_0\}$ of the cluster earthquake map $E(g_0,-)$ at $g_0$ is a covering map. Since the degree of $E(g_0,-)^*$ varies continuously on $g_0$, the degree is a constant. We will prove that this degree is one.

If this degree is more than one, then there is $v\neq v_0\in\bExch_\bs$ such that $\interior \cC^+_{(v)} \cap \interior \cC^+_{(v_0)} = \emptyset$ but that $\interior \cD^+_{(v)}\cap \interior \cD^+_{(v_0)}\neq \emptyset$. Therefore, all we have to do is $\interior \cD^+_{(v)}\cap \interior \cD^+_{(v_0)}= \emptyset$ for any $v\neq v_0\in\bExch_\bs$.


We fix $v_0 \in \bExch_\bs$ and choose a skew-symmetrizable matrix $\varepsilon^{(v_0)}$ satisfying \eqref{eq:Dynkin_quiver} for the Cartan matrix in each case (\cref{cones}).
We consider the point $g_0$ such that $(X^{(v_0)}_1(g_0),X^{(v_0)}_2(g_0))=(1,1)$.
In each case, we define the vertices $\{v_k\}_{k \geq 0}$ of $\bExch_\bs$ so that $v_{2i} \overbar{1} v_{2i+1}$ and $v_{2i+1} \overbar{2} v_{2i+2}$ for $i \geq 0$.

For $A_1\times A_1$, it is easy to see that each $\log \mathbf{X}^{(v_0)}(\cD^+_{(v)})$ corresponds to the first through fourth quadrants, hence $\interior \cD^+_{(v)}(g_0)\cap\interior \cD^+_{(v_0)}(g_0)= \emptyset$ for any $v\neq v_0\in\bExch_\bs$ such that $\interior \cC^+_{(v)} \cap \interior \cC^+_{(v_0)} = \emptyset$.

For $A_2$, we have the vertices  $v_0\overbar{1} v_1\overbar{2} v_2 \overbar{1} \cdots \overbar{1}v_5 \overbar{(1\ 2)} v_0$. 
By \eqref{def_D}, we get
\begin{align*}
    \interior \cD^+_{(v_0)}(g_0)&=\left\{(X_1,X_2)\in\mathbb{R}_{>0}^2\ \middle|\ 
    X_1>1 \text{ and }X_2>1
    \right\},\\
    \interior \cD^+_{(v_1)}(g_0)&=\left\{(X_1,X_2)\in\mathbb{R}_{>0}^2\ \middle|\  \frac{1}{X_1}>1 \text{ and } \frac{X_1\,X_2}{X_1+1}>\frac{1}{2}
    \right\},\\
    \interior \cD^+_{(v_2)}(g_0)&=\left\{(X_1,X_2)\in\mathbb{R}_{>0}^2\ \middle|\ \frac{X_2}{X_1\,X_2+X_1+1}>\frac{1}{3} \text{ and } \frac{X_1+1}{X_1\,X_2}>2
    \right\},\\
    \interior \cD^+_{(v_3)}(g_0)&=\left\{(X_1,X_2)\in\mathbb{R}_{>0}^2\ \middle|\ \frac{X_1\,X_2+X_1+1}{X_2}>3 \text{ and } \frac{1}{X_1\,(X_2+1)}>\frac{1}{2}
    \right\},\\
    \interior \cD^+_{(v_4)}(g_0)&=\left\{(X_1,X_2)\in\mathbb{R}_{>0}^2\ \middle|\ \frac{1}{X_2}>1 \text{ and } X_1\,(X_2+1)>2
    \right\}.
\end{align*}
We can calculate that $X_1<1$ on $\interior \cD^+_{(v_1)}(g_0)\cup\interior \cD^+_{(v_2)}(g_0)$ and that $X_2<1$ on $\interior \cD^+_{(v_3)}(g_0)\cup\interior \cD^+_{(v_4)}(g_0)$. Therefore, we prove that $\interior \cD^+_{(v)}(g_0)\cap\interior \cD^+_{(v_0)}(g_0)= \emptyset$ for any $v\neq v_0\in\bExch_\bs$ such that $\interior \cC^+_{(v)} \cap \interior \cC^+_{(v_0)} = \emptyset$.

For $B_2$ and $G_2$, their proofs are similar. (See also \cref{fig:rk2_eq_plots}.)
\end{proof}

We are going to prepare some notations for the \emph{cluster reduction} argument for earthquake maps. 
Let us choose a vertex $v_0 \in \bExch_\bs$ and a subset $J \subset I$. It corresponds to a face $F:= \{L \in \cC^+_{(v_0)} \mid x_j^{(v_0)}(L)=0 \mbox{ for $j \in J$}\}$. Associated to such a face $F$ is the mutation class $\bs_F$ which contains the seed $((\ve_{ij}^{(v_0)})_{i,j \in J},
(X_i^{(v_0)})_{i \in J})$. The complete subgraph of $\bExch_\bs$ containing the vertices $v$ such that $F \subset \cC^+_{(v)}$ is naturally identified with $\bExch_{\bs_F}$ by \cite[Theorem 5.8]{FG09}\footnote{This is proved for any mutation class by Cao--Li \cite{CL}.}. Therefore we regard $\bExch_{\bs_F} \subset \bExch_\bs$.

Then we have the \emph{cluster reduction morphisms}
\begin{align*}
    \pi_F^+&: \X_\bs(\pos) \to \X_{\bs_F}(\pos), \\
    \pi_F^\trop&: \X_\bs(\bR^\trop) \to \X_{\bs_F}(\bR^\trop)
\end{align*}
obtained by the natural projections
$(X_i^{(v)})_{i \in I} \mapsto (X_i^{(v)})_{i \in J}$ and $(x_i^{(v)})_{i \in I} \mapsto (x_i^{(v)})_{i \in J}$ for $v \in \bExch_{\bs_F}$, respectively.

Let $D_F:=\bigcup_{v \in \bExch_{\bs_F}} \cC^+_{(v)} \subset \X_\bs(\bR^\trop)$ denote the union of the maximal dimensional cones which contain $F$ (the \emph{star neighborhood} of $F$). Let $\Sigma_F$ be the union of the faces $F^\bot_{(v)}:=\{L \in \cC^+_{(v)} \mid x_i^{(v)}(L)=0 \mbox{ for $i \in I \setminus J$}\}$ transverse to $F$ for all $v \in \bExch_{\bs_F}$ (the \emph{link} of $F$). Then $\Sigma_F$ is a PL hypersurface in $\X_\bs(\bR^\trop)$.

\begin{lem}[cluster reduction of the cluster earthquake maps]\label{lem:eq_reduction}
For any face $F$ as above, we have the following commutative diagram:
\begin{equation}\label{eq:reduction_diagram}
    \begin{tikzcd}[column sep=large]
    D_F \ar[r,"{E(g_0,-)|_{D_F}}"] \ar[d,"\pi_F^\trop|_{D_F}"'] & \X_\bs(\pos) \ar[d,"\pi_F^+"] \\
    \X_{\bs_F}(\bR^\trop) \ar[r,"{E(\pi_F^+(g_0),-)}"'] & \X_{\bs_F}(\pos)
    \end{tikzcd}
\end{equation}
for any basepoint $g_0 \in \X_\bs(\pos)$. Moreover, the hypersurface $\Sigma_F$ defines a PL section of $\pi_F^\trop$. 
\end{lem}

\begin{proof}
Since the cones containing $F$ are related to each other by mutations directed to $J$ \cite[Theorem 5.8]{FG09}, the first statement follows from the simple observation that the cluster transformation formulas of the coordinates $X_j^{(v)}$ for $j \in J$ associated with such mutations do not involve the coordinates $X_i^{(v)}$ for $i \in I \setminus J$. 

Since $\Sigma_F$ is the union of the faces $F^\bot_{(v)}$ for all $v \in \bExch_{\bs_F}$, the restriction $\pi_F^\trop|_{\Sigma_F}$ is a PL isomorphism onto its image. Indeed, one can see that the image of each face $F^\bot_{(v)}$ is exactly the cone $\cC^+_{(v)}$ in the fan $\fF^+_{\bs_F}$ again by the observation on cluster transformations. 
Thus the second assertion is proved. 
\end{proof}
Notice that the mutation class $\bs_F$ is of rank $(N-\dim F)$.

\begin{proof}[Proof of \cref{thm:Eq_finite_type}]
We give a proof by induction on the rank $N$ of the mutation class. From \cref{lem:homeo_rank2}, the basis steps $N=1,2$ are established. 

Therefore we consider the case $N \geq 3$, and suppose that the cluster earthquake maps are homeomorphisms for any basepoint $g_0$ and any rank $< N$ mutation classes of finite type. We first claim that the cluster earthquake maps are local homeomorphisms for all rank $N$ mutation classes of finite type. Let $\bs$ be a mutation class of rank $N$ of finite type, and choose a basepoint $g_0 \in \X_\bs(\pos)$. Then $E(g_0,-)$ is obviously a local homeomorphism at any interior point of a maximal dimensional cone in $\fF^+_\bs$. Therefore we need to verify that $E(g_0,-)$ is a local homeomorphism at each point on a codimension $k$ face $F$ in $\fF^+_\bs$. 

Consider the mutation class $\bs_F$ as in \cref{lem:eq_reduction}. Assume first that $F$ is not a point, so that $\bs_F$ has rank $k < N$. Then the earthquake map $E(\pi_F^+(g_0),-)$ is a homeomorphism by the induction assumption. By  \cref{lem:eq_reduction}, 
\begin{align}\label{eq:PL_bundle}
    \pi_F^\trop|_{D_F}: D_F \cong \Sigma_F \times F \to \X_{\bs_F}(\bR^\trop)
\end{align}
is a PL bundle with fiber $F$. Indeed, notice first that if we fix a vertex $v_0 \in \bExch_{\bs_F}$, for $i \in I \setminus J$, the coordinates $x_i^{(v)}$  associated to $v \in \bExch_{\bs_F}$ can be written as $x_i^{(v)}=x_i+f_i((x_j)_{j \in J})$ for some linear functions $f_i$, where we write $x_i:=x_i^{(v_0)}$, again by the observation on the cluster $\X$-transformation. 
Hence the fiber over a point $\xi=(\xi_j)_{j \in J} \in \cC^+_{(v)} \subset \X_{\bs_F}(\bR^\trop)$ is expressed as
\begin{align*}
    (\pi_F^\trop|_{D_F})^{-1}(\xi)=\{ (x_i)_{i \in I} \mid x_j=\xi_j~\mbox{for $j \in J$},~x_i+f_i(\xi)\geq 0~\mbox{for $i \in I \setminus J$}\},
\end{align*}
which is identified with $F$ via the linear map 
\begin{align*}
    (\pi_F^\trop|_{D_F})^{-1}(\xi) \xrightarrow{\sim} F, \quad ((x_i)_{i \in I\setminus J},\, (\xi_j)_{j \in J}) \mapsto ((x_i+f_i(\xi))_{i \in I\setminus J}, (0)_{j \in J}).
\end{align*}
These identifications combine to give a PL local trivialization of \eqref{eq:PL_bundle}. 
Moreover, $E(g_0,-)$ maps distinct fibers of \eqref{eq:PL_bundle} to disjoint fibers of $\pi_F^+$, since the base map $E(\pi_F^+(g_0),-)$ is injective. 
Indeed, if there is an intersection point $p$ of the images of fibers over two points $\xi_1, \xi_2 \in \cX_{\bs_F}(\bR^\trop)$, then $\xi_1$ and $\xi_2$ must be send to the same point $\pi^+_F(E(g_0, p))$ by $E(\pi^+_F(g_0), -)$, which contradicts to the injectivity.
Note also that each fiber $(\pi_F^\trop|_{D_F})^{-1}(\xi)$ is contained in $\cC^{+}_{(v)} \subset \X_\bs(\bR^\trop)$, and hence the restriction of $E(g_0,-)$ to each fiber is injective. 
Therefore $E(g_0,-)|_{D_F}$ is a local homeomorphism at any point on the central fiber $F$ by the Brouwer's invariance of domain theorem. Since $F$ is an arbitrary face of dimension $>0$, it follows that 
\begin{align*}
    E(g_0,-)^\ast:=E(g_0,-)|_{\X_\bs(\bR^\trop) \setminus \{0\}}: \X_\bs(\bR^\trop) \setminus \{0\} \to \X_\bs(\pos)\setminus \{g_0\}
\end{align*}
is a proper local homeomorphism. Since the domain and the target spaces are Hausdorff and locally compact (as they are topologically Euclidean spaces minus one point), it implies that $E(g_0,-)^\ast$ is a covering map. It must be a homeomorphism, since $\X_\bs(\pos) \setminus \{g_0\} \simeq S^{N-1}$ is simply-connected for $N \geq 3$. Filling the missing point, we see that $E(g_0,-)$ is a homeomorphism. 
Thus the assertion is proved. 
\end{proof}

\section{Derivatives of cluster earthquake maps}
In this section, we study the derivatives of the earthquake maps for a seed of finite type.

\subsection{Tangent bundles of cluster manifolds}\label{sec:tangent}

We are going to construct an isomorphism between the tangent bundle $T \cX_\bs(\bR_{>0})$ of the cluster manifold  and the trivial topological fiber bundle $\cX_\bs(\bR_{>0}) \times \cX_\bs(\bR^\trop)$ over $\cX_\bs(\bR_{>0})$ via a map similar to \eqref{eq:quakeflow}. 

\begin{thm}\label{thm:tangent_trop}
For any mutation class $\bs$ of finite type, the tangent map
\begin{align}\label{eq:tangent}
    dE: \cX_\bs(\bR_{>0}) \times \cX_\bs(\bR^\trop) \to T\cX_\bs(\bR_{>0})\, ,\ 
    (g, L) \mapsto \Big(g,\ \frac{d}{dt}\bigg|_{t=0^+} E(g,\, t \cdot L) \Big)
\end{align}
is a $\Gamma_\bs$-equivariant isomorhpism of topological fiber bundles.
Here, $t\cdot$ denotes the natural $\bR_{>0}$-action on $\cX_\bs(\bR^\trop)$. 
\end{thm}

Here, note that the natural PL structure of $\X_\bs(\bR^\trop)$ induces a \emph{tangential structure} in the sense of \cite{Bo98}. The one-sided directional derivative in the statement is well-defined in this sense. See \cite[Section 1]{Bo98} for a detail.

We will show \cref{lem:tangent_rank2} for the induction for the other proof of \cref{thm:tangent_trop} than  $\Gamma_\bs$-equivalence.

\begin{lem}\label{lem:tangent_rank2}
For any rank two mutation classes of finite type, the tangent map \eqref{eq:tangent}
is an isomorhpism of topological fiber bundles.
\end{lem}

\begin{proof}
By \cref{thm:finite_seed_pattern}, it suffices to consider the cases of $A_1 \times A_1$, $A_2$, $B_2$ and $G_2$. We fix a point $g\in\cX_\bs(\bR_{>0})$ and let $dE_g:\cX_\bs(\bR^{\trop})\to T_g\cX_\bs(\bR_{>0})$ be the restriction of $dE$ at $g$. We will prove that $dE_g$ is a homeomorphism.

For  $L\in\cC_{(v)}^+$, by using the coordinate $\log \mathbf{X}^{(v)}:=\left(\log X^{(v)}_1,\ \log X^{(v)}_2 \right)$, 
\begin{align*}
    \frac{d}{dt}\bigg|_{t=0^+} \log \textbf{X}^{(v)}( E(g,\, t \cdot L))
    &=x_1^{(v)}(L) \bigg(\frac{\partial}{\partial \log X_1^{(v)}}\bigg)_{\!\!g} +x_2^{(v)}(L) \bigg(\frac{\partial}{\partial \log X_2^{(v)}}\bigg)_{\!\!g}.
\end{align*}
Therefore, $dE_g|_{\cC_{(v)}^+}$ is a homeomorphism from $\{L\in\cX_\bs(\bR^{\trop})\mid x_1^{(v)}(L)\geq0,\ x_2^{(v)}(L)\geq0\}$ to 
$\bigg\{a_1\bigg(\dfrac{\partial}{\partial \log X_1^{(v)}}\bigg)_{\!\!g} + a_2\bigg(\dfrac{\partial}{\partial \log X_2^{(v)}}\bigg)_{\!\!g} \in T_g\cX_\bs(\bR_{>0}) \ \Big|\ a_1\geq0,\ a_2\geq0\bigg\}$. 


For the type $A_2$, we have the vertices  $v_0\overbar{1} v_1\overbar{2} v_2 \overbar{1} \cdots \overbar{1}v_5 \overbar{(1\ 2)} v_0$.
We choose some point $l_k$ from $\cC^+_{(v_k)} \cap \cC^+_{(v_{k+1})}$ for $k \in \bZ/5\bZ$.
For example, we choose $l_0,...,l_4$ to be $(1,0),\ (0,1),\ (-1,0),\ (0,-1), (1,-1)$ with respect to the coordinate $\mathbf{x}^{(v_0)} = (x^{(v_0)}_1, x^{(v_0)}_2)$, respectively in this order.
We choose $g$ such that $X_1^{(v_0)}(g)=1$ and $X_2^{(v_0)}(g)=1$, and define $\xi_k:=dE_g(l_k)$.
Then, we directly compute $\xi_0,...,\xi_4=(1,0),\ (0,1),\ (-1,1/2),\ $
$(-2/3,-2/3),\ (1/2, -1)$ with respect to $\bigg(\dfrac{\partial}{\partial \log X_1^{(v)}}\bigg)_{\!\!g}$ and $\bigg(\dfrac{\partial}{\partial \log X_2^{(v)}}\bigg)_{\!\!g}$, respectively. 
The angles of $l_k$ determine the cyclic order of $(l_k)_{k=0}^{4}$. The same applies to $(\xi_k)_{k=0}^{4}$. 
We see that $dE_g$ preserves the cyclic orders of $(l_k)_{k=0}^{4}$ and $(\xi_k)_{k=0}^{4}$. Additionally, each restriction $dE_g|_{\cC_{(v_i)}^+}$ is a homeomorphism from each cone between adjacent vectors of $(l_k)_{k=0}^{4}$  to  each cone between adjacent vectors of $(\xi_k)_{k=0}^{4}$. Therefore, $dE_g$ is a homeomorphism. When $X_1^{(v_0)}(g)\neq1$ or $X_2^{(v_0)}(g)\neq1$, $dE_g$ preserves the cyclic order similarly and hence is a homeomorphism for the type $A_2$.

Below we show a table, including the other rank two cases of finite type, from which we can see that $dE_g$ preserves the cyclic order.
Here, the vectors $l_k$ are chosen to be in $\cC_{(v_k)}^+ \cap \cC^+_{(v_{k+1})}$ and to be easy to calculate, and we define $\xi_k: = dE_g(l_k)$, as above.

\vspace{7pt}
\setlength{\extrarowheight}{10pt}
\begin{spacing}{0.5}
\begin{tabular}{|c|c|c|}
\hline
&$l_k$ & $\xi_k$  \\
\hline
$A_1\times A_1$ &$(1,0),\ (0,1),\ (-1,0),\ (0,-1)$ &$ (1,0),\ (0,1),\ (-1,0),\ (0,-1)$\\
\hline
$A_2$ &\begin{tabular}{c}
        $(1,0),\ (0,1),\ (-1,0),\ (0,-1),$\\
        $(1,-1)$
       \end{tabular}
      &\begin{tabular}{c}
          $(1,0),\ (0,1)\ (-1,\frac{1}{2}),\ (-\frac{2}{3},-\frac{2}{3}),$\\
          $(\frac{1}{2},-1)$
      \end{tabular}\\
\hline
$B_2$ &\begin{tabular}{c}
            $(1,0),\ (0,1),\ (-1,0),\ (0,-1),$\\
            $(1,-2),\ (1, -1)$ 
        \end{tabular}
      &\begin{tabular}{c}
        $(1,0),\ (0,1),\ (-1,1),\ (-\frac{4}{5}, -\frac{1}{5}),$\\
        $(-\frac{1}{3},-\frac{4}{3}),\ (\frac{1}{2},-1)$
       \end{tabular}\\
\hline
$G_2$ &\begin{tabular}{c}
        $(1,0),\ (0,1),\ (-1,0),\ (0,-1),$\\
        $(1,-3),\ (1,-2),\ (2,-3),\ (1,-1)$
       \end{tabular}
      &\begin{tabular}{c}
        $(1,0),\ (0,1),\ (-1,\frac{3}{2}),\ (-\frac{8}{9},\frac{1}{3}),$\\ 
        $(-\frac{7}{5},-\frac{3}{5}),\ (-\frac{1}{2},-\frac{13}{14}),\ (0,-2),\ (\frac{1}{2},-1)$ 
       \end{tabular}\\
\hline
\end{tabular}
\end{spacing}
\vspace{20pt}
Therefore, $dE_g\ (g\in\cX_\bs(\bR_{>0}))$ is a homeomorphism and $dE$ is an isomorphism of topological fiber bundles for any rank two mutation classes of finite type.
\end{proof}

\begin{proof}[Proof of \cref{thm:tangent_trop}]
    
We show that the tangent map
\begin{align*}
    dE(g,-)_0: T_0\cX_\bs(\bR^\trop) \to T_g \cX_\bs(\bR_{>0})
\end{align*}
is a homeomorphism for any $g \in \cX_\bs(\bR_{>0})$. Here, $T_0\cX_\bs(\bR^\trop)$ is the tangent space at the origin of the tangentiable manifold $\cX_\bs(\bR^\trop)$, namely, the set of one-sided directional derivatives (\cite[Section 1]{Bo98}).
Let us proceed by induction on the rank $N \geq 1$ of the mutation class $\bs$, producing a \emph{tangential version} of the cluster reduction argument in the proof of \cref{thm:Eq_finite_type}. The rank one case is obvious, and the rank two cases are verified in \cref{lem:tangent_rank2}. We consider the case $N \geq 3$, and assume that the statement holds true for any basepoint $g$ and any rank $<N$ mutation classes of finite type. Take the tangent maps of the diagram \eqref{eq:reduction_diagram}:
\begin{equation*}
    \begin{tikzcd}[column sep=large]
    T_0D_F \ar[r,"{dE(g,-)_0}"] \ar[d,"{(d\pi_F^\trop)_0}"'] & T_g\X_\bs(\pos) \ar[d,"{(d\pi_F^+)_g}"] \\
    T_0\X_{\bs_F}(\bR^\trop) \ar[r,"{dE(g',-)_0}"'] & T_{g'}\X_{\bs_F}(\pos),
    \end{tikzcd}
\end{equation*}
where $g':=\pi_F^+(g)$. It is also commutative by the functoriality of tangent maps. Observe that via the canonical identifications $T_0 D_F \cong D_F$ and $T_0 \X_{\bs_F}(\bR^\trop) \cong \X_{\bs_F}(\bR^\trop)$, the tangent map $(d\pi_F^\trop)_0$ is exactly the same as $\pi_F^\trop$ (just similarly to the relation between a linear map and its tangent map). Hence $(d\pi_F^\trop)_0$ is also a fiber bundle with fiber $F$. Therefore we can apply the same argument as in the proof of \cref{thm:Eq_finite_type}, concluding that the tangent map $dE(g,-)_0$ is a local homeomorphism at any point on the central fiber $F$. It follows that  $dE(g,-)_0^\ast: T_0\cX_\bs(\bR^\trop) \setminus\{0\} \to T_g \cX_\bs(\bR_{>0})\setminus\{0\}$ is a covering map, so must be a homeomorphism. Thus the assertion is proved. 

Also, the $\Gamma_\bs$-equivariance is verified as
\begin{align*}
    \phi(dE(g,L)) &= \Big(\phi(g),\ \frac{d}{dt}\bigg|_{t=0^+}\phi(E(g,\, t \cdot L)) \Big) \\
    &= \Big(\phi(g),\ \frac{d}{dt}\bigg|_{t=0^+}E(\phi(g),\, t \cdot \phi(L)) \Big) \\
    &=dE(\phi(g),\phi(L)).
\end{align*}

\end{proof}

\subsection{The fans in the tangent spaces of cluster manifolds}
Let us fix a mutation class $\bs$
of finite type 
and a vertex $v_0 \in \bExch_\bs$.
We recall that the cones $x_{(v_0)}(\cC^+_{(v)})$ for $v \in \bExch_\bs$ gives a complete fan $x_{(v_0)}(\fF^+_{\bs})$ in $\bR^I$.

For $g \in \cX_\bs(\bR_{>0})$, the tangent map 
\begin{align*}
    dE(g,-)_0: T_0 \cX_\bs(\bR^\trop) \cong \cX_\bs(\bR^\trop) \to T_{g} \cX_\bs(\bR_{>0})
\end{align*}
induces a fan $\fF^+_{\bs}(g)$ in $T_{g}\cX_\bs(\pos)$
given by the faces of the cones
\begin{align*}
    \cC^+_{(v)}(g) := dE(g,-)_0(\cC^+_{(v)}) \subset T_{g}\cX_\bs(\bR_{>0})
\end{align*}
for $v \in \bExch_\bs$.

Since the map $dE(g,-)_0$ is bijective as we saw in the proof of \cref{thm:tangent_trop}, $(d \log \mathbf{X}^{(v_0)})_{g} (\fF^+_{\bs}(g))$ is a complete fan in $\bR^I$.
The following theorem states that although this fan is not quite the same as the fan $x_{(v_0)}(\fF^+_{\bs})$,
it is asymptotically close to the latter as $g$ diverges toward $\interior \bS \cC^+_{(v_0)}$.

\begin{thm}\label{thm:fan_in_tangent}
Let $\bs$ be a mutation class of finite type, and fix $v_0 \in \bExch_\bs$.
Then, 
\begin{align*}
    (d \log \mathbf{X}^{(v_0)})_{g}(\fF^+_{\bs}(g)) \to x_{(v_0)}(\fF^+_{\bs}) \quad \mbox{as $g$ diverges toward $\interior \bS \cC^+_{(v_0)}$}.
\end{align*}
Namely, the cones $(d \log \mathbf{X}^{(v_0)})_{g} (\cC^+_{(v)}(g))$ converges to $x_{(v_0)}(\cC^\pm_{(v)})$ in the Gromov--Hausdorff sense as $g \in \cX_\bs(\pos)$ diverges so that $\mathbf{X}^{(v_0)}(g)^{-1} = (X^{(v_0)}_i (g)^{-1})_{i \in I} \to \mathbf{0}$.
\end{thm}

We postpone the proof of this theorem until the next section.

\section{Asymptotic behavior of the cluster earthquake maps}\label{sec:asymptotics}

For a point $g_0 \in \cX_\bs(\pos)$, the set of intersections of the subsets $\cD^+_{(v)}(g_0)$ for $v \in \bExch_\bs$ is a simplicial set but not a fan in the log-coordinate chart $\log \mathbf{X}^{(v_0)}(\cX_\bs(\pos))$ for any $v_0 \in \bExch_\bs$. 
Namely, the simplices are not cones, as they are not invariant under the $\pos$-action.
(Roughly speaking, the boundary of simplices are out of linear.)
Nevertheless, we are going to see that the set is asymptotically a fan at the \emph{Thurston boundary} of $\cX_\bs(\pos)$ in a certain sense.

Recall that the subsets $\cD^+_{(v)}(g_0)$ consist of the images $E(g_0, L)$ of the points $g_0 \in \cX_\bs(\pos)$ and $L \in \cX_\bs(\bR^\trop)$ under the cluster earthquake map $E: \cX_\bs(\pos) \times \cX_\bs(\bR^\trop) \to \cX_\bs(\pos)$.
Letting $L$ go to infinity, 
this point goes to the Thurston boundary of $\cX_\bs(\pos)$.
We will investigate this direction of limit in \cref{subsec:L_to_infty}.
On the other hand, one can also let $g_0$ go to infinity. 
We will investigate this direction of limit in \cref{subsec:g_to_infty}.

\subsection{The case where $L$ goes to infinity}\label{subsec:L_to_infty}
\begin{figure}[h]
    \centering
    \begin{tikzpicture}
    \node at (0,0) {\includegraphics[width=4cm]{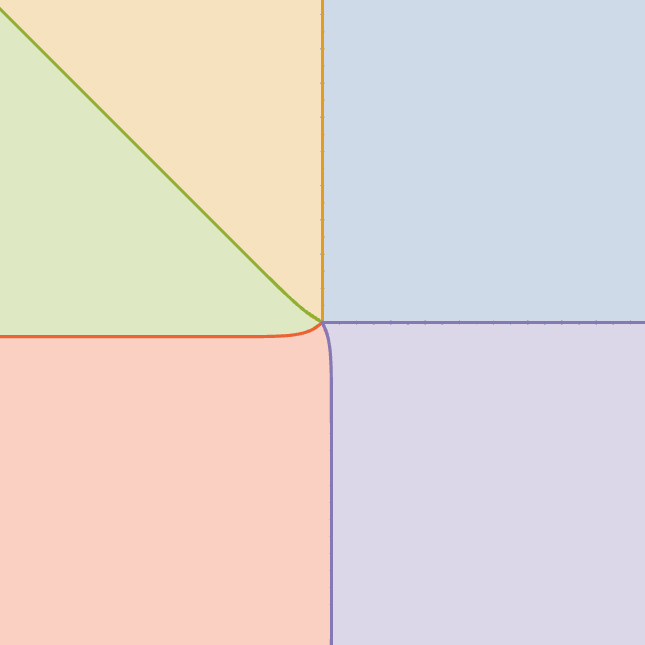}};
    \node at (5,0) {\includegraphics[width=4cm]{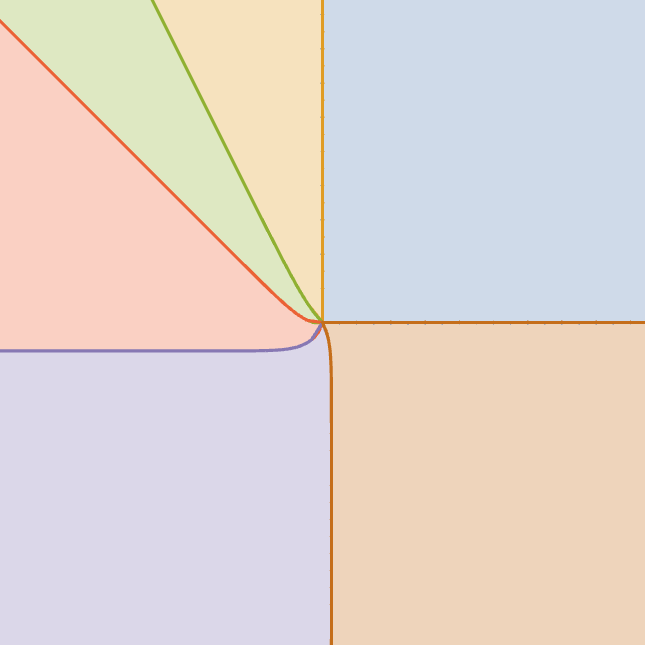}};
    \node at (10, 0) {\includegraphics[width=4cm]{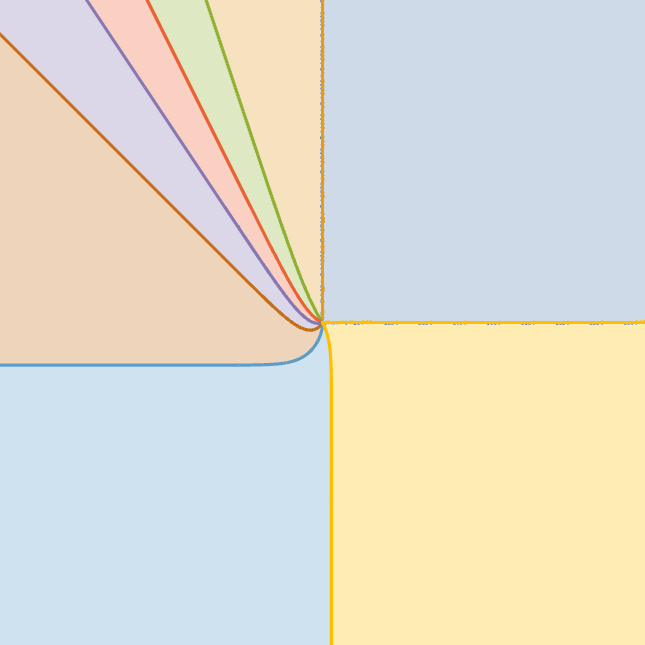}};
    \node at (0,-2.5) {$A_2$};
    \node at (5,-2.5) {$B_2$};
    \node at (10,-2.5) {$G_2$};
    \end{tikzpicture}
    \caption{The simplices $\log \mathbf{X}^{(v_0)} (E(g^{(v_0)}_1, \cC^+_{(v)}))$ for $v \in \bExch_\bs$. These plots are drawn in the wider range $[-25, 25] \times [-25, 25]$ than the plots in the proof of \cref{lem:homeo_rank2}.}
    \label{fig:log_earthquake_far}
\end{figure}

\begin{figure}[h]
    \centering
    \begin{tikzpicture}
    \node at (0,0) {\includegraphics[width=4cm]{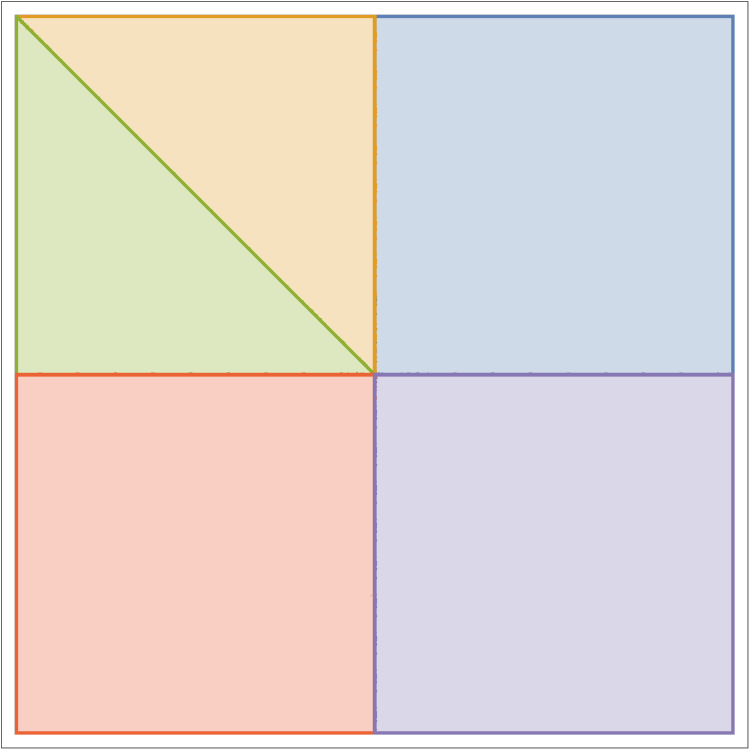}};
    \node at (5,0) {\includegraphics[width=4cm]{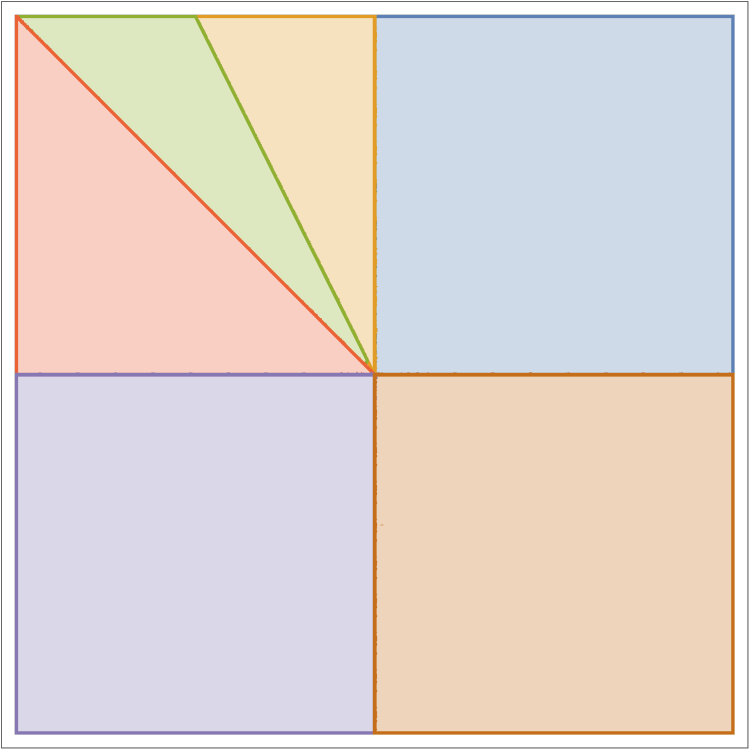}};
    \node at (10, 0) {\includegraphics[width=4cm]{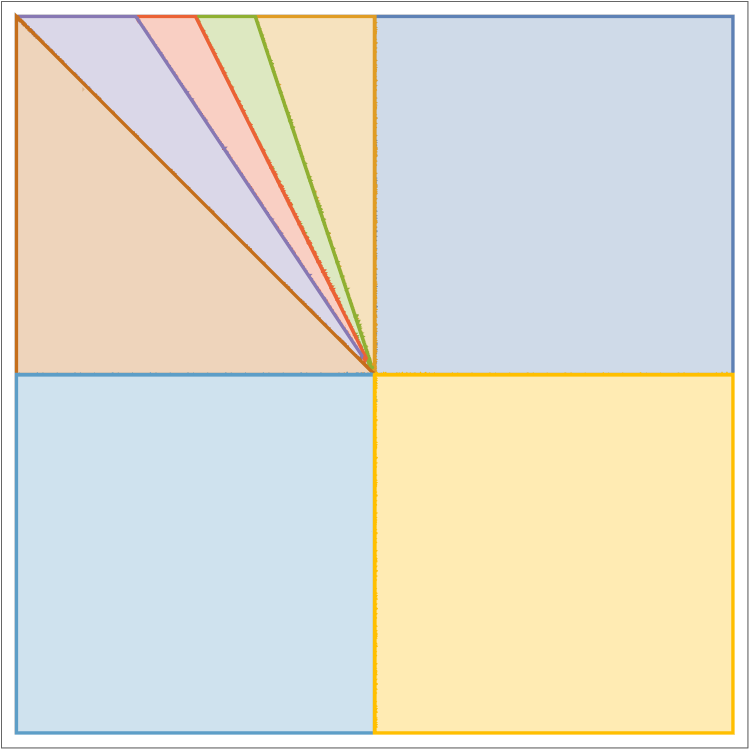}};
    \node at (0,-2.5) {$A_2$};
    \node at (5,-2.5) {$B_2$};
    \node at (10,-2.5) {$G_2$};
    \end{tikzpicture}
    \caption{The Fock--Goncharov fan of the opposite mutation class $-\bs$.
    Namely, the cones $\mathbf{x}^{(v_0)}(\iota^{-1}(\cC^+_{(-v)}))$ in $\bR^I$.}
    \label{fig:cluster_complex_minus}
\end{figure}

Fix a mutation class $\bs$ of finite type, a point $g_0 \in \cX_\bs(\pos)$ and a vertex $v_0 \in \bExch_\bs$.
Let us consider the plot $\log \mathbf{X}^{(v_0)} \circ E(g_0,-):\X_\bs(\bR^\trop) \to \bR^I$ of the cluster earthquake map in the log-coordinates associated with $v_0$. 
Our preliminary observation is that the plot shown in \cref{fig:log_earthquake_far} seen from a distance looks similar to the Fock--Goncharov fan of the opposite mutation class $-\bs$ (\cref{def:opposite}), as shown in \cref{fig:cluster_complex_minus}. Let us make it more precise.


%
%

There is a PL isomorphism 
\begin{align}\label{eq:iota}
    \iota: \cX_\bs(\bR^\trop) \xrightarrow{\sim} \cX_{-\bs}(\bR^\trop)
\end{align}
such that $x_i^{(-v)}(\iota(L)) = -x_i^{(v)}(L)$
for each $v \in \bExch_\bs$ and $L \in \cX_\bs(\bR^\trop)$
\footnote{One can verify that there is an isomorphism $\iota: \cX_\bs(\pos) \xrightarrow{\sim} \cX_{-\bs}(\pos)$ in the same way.}.
We will identify $\cX_\bs(\bR^\trop)$ and $\cX_{-\bs}(\bR^\trop)$ via this isomorphism $\iota$.
Let us recall the \emph{Thurston compactification} $\overline{\cX_\bs(\pos)}$ of a cluster manifold $\cX_\bs(\pos)$ \cite{Le16,FG16,Ish19}. It is given by the spherical compactification
\begin{align*}
    \overline{\cX_\bs(\pos)} := \cX_\bs(\pos) \sqcup \bS \cX_\bs(\bR^\trop)
\end{align*}
where $\bS \cX_\bs(\bR^\trop) := (\cX_\bs(\bR^\trop) \setminus \{0\}) / \pos$.
It is endowed with a topology such that a divergent sequence $g_n$ in $\cX_\bs(\pos)$ converges to $[L] \in \bS \cX_\bs(\bR^\trop)$ in $\overline{\cX_\bs(\pos)}$ if and only if $[\log \mathbf{X}^{(v)}(g_n)]$ converges to $[\mathbf{x}^{(v)}(L)]$ in $\bS \bR^I=(\bR^I \setminus \{0\})/\bR_{>0}$  (cf. \cite[Section 7]{Le16}).

For $g_0 \in \cX_\bs(\pos)$ and $v \in \bExch_\bs$, we define
\begin{align*}
    \cE^+_{(v)}(g_0) := \Big\{ \lim_{t \to \infty} E(g_0,\, t \cdot L)\ \Big|\ L \in \cC^+_{(v)} \Big\} \subset \bS \cX_\bs(\bR^\trop).
\end{align*}

\begin{thm}\label{thm:D_at_bdry}
Under the identification \eqref{eq:iota}, we have
\begin{align*}
    \cE^+_{(v)}(g_0) = \iota^{-1}(\bS \cC^+_{(-v)}) \subset \bS \cX_\bs(\bR^\trop)
\end{align*}
for each $g_0 \in \cX_\bs(\pos)$ and $v \in \bExch_\bs$.
Thus, the set of the faces of $\{\cE^+_{(v)}(g_0) \mid v \in \bExch_\bs\}$ coincides with the cluster complex \cite{FZ-CA2}
$\bS \fF^+_{-\bs} := \{ \bS \cC \mid \cC \in \fF^+_{-\bs}\}$ via the identification $\iota$.
\end{thm}
\begin{proof}
Take $v \in \bExch_\bs$.
For $k \in I$, let $L = L^{(v)}_{k} \in \cC_{(v)}^+$ be the elementary lamination characterized by $x_i^{(v)}(L^{(v)}_k) = \delta_{ik}$. Then from the definition, the earthquake along $t \cdot L$ satisfies $X_i^{(v)}(E(g_0,t\cdot L)) = e^{\delta_{ik}t} X_i^{(v)}(g_0)$. 
By using the separation formula (\cref{thm:sep_for}), we rewrite it in the coordinates $X_i^{(v_0)}$ as
\begin{align*}
    X_i^{(v_0)}(E(g_0,\, t \cdot L)) 
    &= \Big( \prod_{j \in I} (X_j^{(v)}(E(g_0,\, t \cdot L)))^{c_{ij}} \cdot F_j^{v \to v_0}(\mathbf{X}^{(v)}(E(g_0,\, t \cdot L)))^{\varepsilon^{(v)}_{ij}} \Big)\\
    &= (e^t)^{c_{ik}} \cdot \Big( \prod_{j \in I} (X_j^{(v)}(g_0))^{c_{ij}} \cdot F_j^{v \to v_0}(\mathbf{X}^{(v)}(E(g_0,\, t \cdot L)))^{\varepsilon^{(v)}_{ij}} \Big)
\end{align*}
where $C^{\bs}_{v \to v_0} = (c_{ij})_{i,j \in I}$.
Hence,
\begin{align*}
    \log X_i^{(v_0)}(E(g_0,\, t \cdot L)) &=
    t \cdot c_{ik} +\sum_{j \in I}c_{ij} \log X_j^{(v)}(g_0)
    + \sum_{j \in I} \varepsilon^{(v)}_{ij} \log F^{v \to v_0}_j(\mathbf{X}^{(v)} (E(g_0,\, t \cdot L))).
\end{align*}
Recall that the $F$-polynomial has the following form:
\begin{align}\label{eq:F_exp}
    F_j^{v \to v_0}(\mathbf{X}^{(v)}) = 1 + \sum_{\alpha=(\alpha_m)} \xi_\alpha\cdot \prod_{m \in I} (X_m^{(v)})^{\alpha_m} + \prod_{m \in I} (X_m^{(v)})^{f_{jm}}
\end{align}
for some constants $\xi_\alpha \in \bZ$,
where $F_{v \to v_0}^\bs=(f_{jm})_{j,m \in I}$ is the $F$-matrix (\cref{c_mat} \eqref{item:F_mutation}), and $\alpha=(\alpha_m)$ runs over certain set of integral vectors satisfying $\alpha_m \leq f_{jm}$ for all $m \in I$. 
Then the last term is estimated as 
\begin{align*}
    &\frac{1}{t}\log F^{v \to v_0}_j(\mathbf{X}^{(v)} (E(g_0,\, t \cdot L))\\
    &\hspace{3cm}= \frac{1}{t}\log F_j^{v \to v_0}(X_1^{(v)}(g_0), \dots, e^t X_k^{(v)}(g_0), \dots, X_N^{(v)}(g_0))\\
    &\hspace{3cm}= \frac{1}{t}\log \Big(1 + \sum_{\alpha} e^{t \cdot \alpha_k}\xi_\alpha\cdot \prod_{m} X_m^{(v)}(g_0)^{\alpha_m} + e^{t \cdot f_{jk}} \prod_{m}X^{(v)}_m(g_0)^{f_{jm}}\Big)\\
    &\hspace{3cm}\xrightarrow[t \to \infty]{} f_{jk},
\end{align*}
since $e^{tf_{jk}}$ is the dominant term inside the logarithm. 
Therefore,
\begin{align*}
    \lim_{t \to \infty} \frac{\log X_i^{(v_0)}(E(g_0,\, t \cdot L))}{t}
    = c_{ik} + \sum_{j \in I} \varepsilon_{ij}^{(v)}f_{jk}.
\end{align*}
By \cite[Theorem 2.8]{FuGy}, we have
\begin{align*}
    C^{\bs}_{v \to v_0} + \varepsilon^{(v)} F^{\bs}_{v \to v_0}
    = C^{-\bs}_{v \to v_0}.
\end{align*}
Thus, $\bS \mathbf{x}^{(v_0)}(\cE_{(v)}^+(g_0))$ is the convex hull in $\bS \mathbf{x}^{(v_0)}(\cX_\bs(\bR^\trop))$ spanned by the points represented by the column vectors of $C_{v \to v_0}^{-\bs} = (G_{v_0 \to v}^\bs)^\tr$.
On the other hand, the cone $x_{(-v_0)}(\cC^+_{(-v)})$ is also spanned by the $g$-vectors from $v_0$ to $v$ at $\bs$ by \cref{prop:cl_cpx=g-vect_fan}.
Thus these two cones coincide with each other.
\end{proof}

\subsection{The case where $g_0$ goes to the boundary at infinity}\label{subsec:g_to_infty}
We continue to fix a seed pattern $\bs$ of finite type and a vertex $v_0 \in \bExch_\bs$.
For $g \in \cX_\bs(\pos)$ and $v_0 \in \bExch_\bs$, let us consider the coordinates
\begin{align*}
    \mathbf{u}^{(v)}_{g}(L):= \log \frac{\mathbf{X}^{(v)}(E(g,L))}{\mathbf{X}^{(v)}(g)}
\end{align*}
and the  simplicial set $\mathbf{u}^{(v_0)}_{g}(\fF^+_\bs)$. Observe that
\begin{align*}
    \mathbf{u}^{(v_0)}_{g}(\cC^+_{(v)}) = \log \frac{\mathbf{X}^{(v_0)}}{\mathbf{X}^{(v_0)}(g)}(\cD^+_{(v)}(g)).
\end{align*}
Although $\mathbf{u}^{(v_0)}_{g}(\fF^+_\bs)$ is not a fan (\emph{i.e.}, the subsets $\mathbf{u}^{(v_0)}_{g}(\cC^+_{(v)})$ are not cones),
we are going to see that it converges to the fan $\mathbf{x}^{(v_0)}(\fF^+_\bs)$ as $g$ diverges toward an appropriate direction.

\begin{thm}\label{thm:g_to_infty}
Let $\bs$ be a mutation class of finite type and $v_0 \in \bExch_\bs$.
Then, 
\begin{align*}
    \mathbf{u}^{(v_0)}_{g}(\fF^+_\bs) \to \mathbf{x}^{(v_0)}(\fF^+_\bs)
    \quad \mbox{as $g$ diverges toward $\interior \bS \cC^+_{(v_0)}$}.
\end{align*}
Namely, the subsets $\mathbf{u}^{(v_0)}_{g}(\cC^+_{(v)})$ converge to $\mathbf{x}^{(v_0)}(\cC^+_{(v)})$ in the Gromov--Hausdorff sense as $g  \in \cX_\bs(\pos)$ diverges so that $\mathbf{X}^{(v_0)}(g)^{-1} = (X^{(v_0)}_i (g)^{-1})_{i \in I} \to \mathbf{0}$.
\end{thm}

\begin{figure}[h]
    \centering
    \begin{tikzpicture}
    \node at (0,0) {\includegraphics[width=4cm]{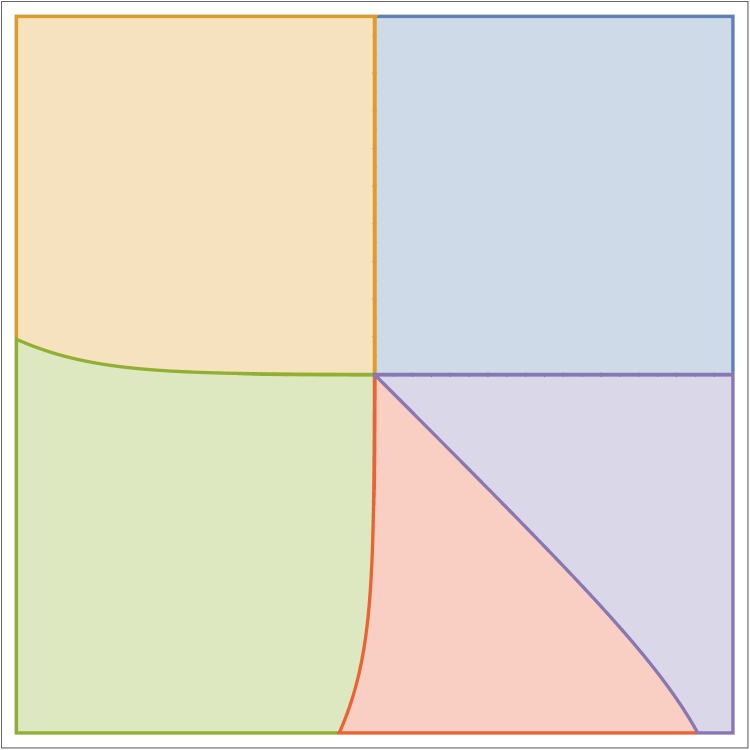}};
    \node at (5,0) {\includegraphics[width=4cm]{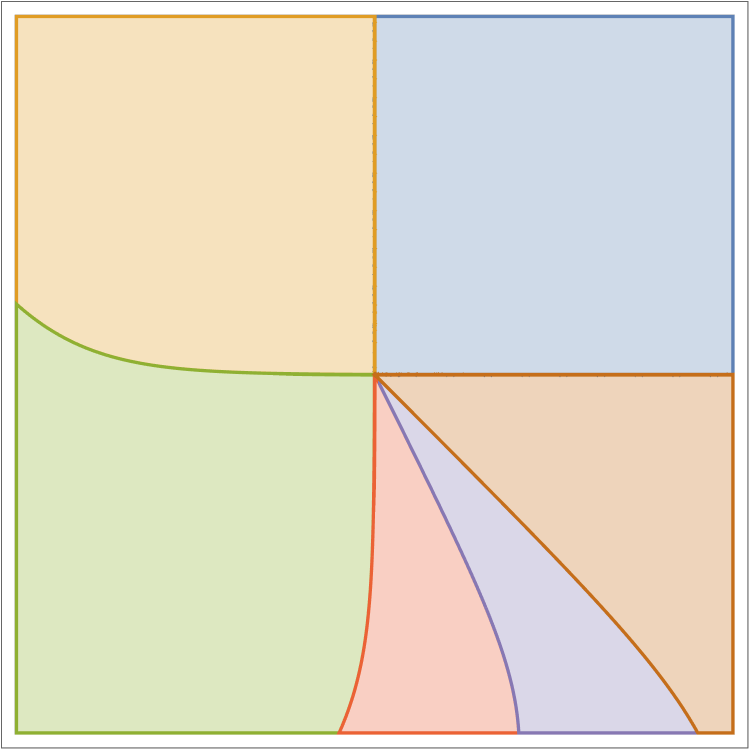}};
    \node at (10, 0) {\includegraphics[width=4cm]{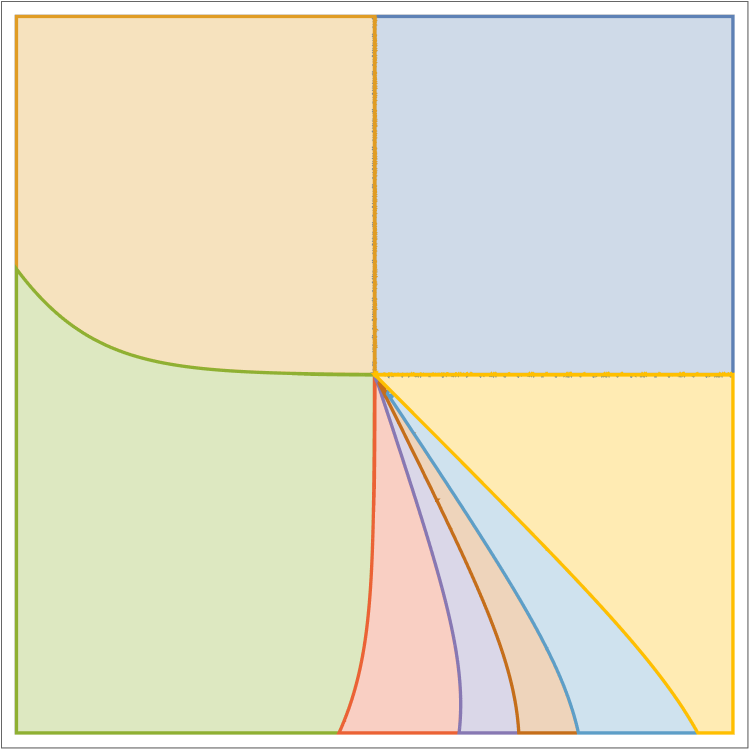}};
    \node at (0,-2.5) {$A_2$};
    \node at (5,-2.5) {$B_2$};
    \node at (10,-2.5) {$G_2$};
    \end{tikzpicture}
    \caption{The plots of the fans $\mathbf{u}_g^{(v_0)}(\fF^+_\bs)$ for $\bs$ are the mutation classes of finite type and rank 2. Here, $\log \mathbf{X}^{(v_0)}(g) = (\log 500, \log 500)$ and they are drawn in the range $[-6, 6] \times [-6, 6]$, which is the same range of the plots in the proof of \cref{lem:homeo_rank2}.}
    \label{fig:u_g(F)}
\end{figure}

Pictures of $\mathbf{u}_g^{(v_0)}(\fF^+_\bs)$ are shown in \cref{fig:u_g(F)}.
Observe how the simplicial set $\mathbf{u}_g^{(v_0)}(\fF^+_\bs)$ converges to $\mathbf{x}^{(v_0)}(\fF^+_\bs)$ by comparing with the plots of $\mathbf{x}^{(v_0)}(\fF^+_\bs)$ in the proof of \cref{lem:homeo_rank2}.

\begin{proof}
For $v \in \bExch_\bs$ and $L \in \cC^+_{(v)}$, we put $T_i := \exp(x^{(v)}_i(L))$.
Let us write
$C^\bs_{v \to v_0} = (c_{ij})_{i,j \in I}$ and $F^\bs_{v \to v_0} = (f_{ij})_{i,j\in I}$.
Then for any $g \in \cX_\bs(\pos)$, by the separation formula (\cref{thm:sep_for}), we compute
\begin{align}
    X_i^{(v_0)}(E(g, L))
    &= \prod_{j \in I} X_j^{(v)}(E(g, L))^{c_{ij}} \cdot F_j^{v \to v_0} (\mathbf{X}^{(v)}(E(g, L)))^{\varepsilon_{ij}^{(v)}}\nonumber\\
    &=\Big(\prod_{j \in I} T_j^{\,c_{ij}}\Big) \cdot \Big(\prod_{j \in I} (X_j^{(v)}(g))^{c_{ij}} \cdot F_j^{v \to v_0}(T \cdot \mathbf{X}^{(v)}(g))^{\varepsilon^{(v)}_{ij}}\Big),\label{eq:sep_with_eq}
\end{align}
where $T \cdot \mathbf{X}^{(v)} = (T_i \cdot X_i^{(v)})_{i \in I}$.
Let $L^{(v)}_i \in \cC^+_{(v)}$ be the point satisfying $x^{(v)}_j(L^{(v)}_i) = \delta_{ij}$ for $i \in I$.
Then, we have
\begin{align*}
    u^{(v_0)}_{g, i}(L^{(v)}_k) = 
    \log \frac{X_i^{(v_0)}(E(g, L^{(v)}_k))}{X_i^{(v_0)}(g)}
    &=
    c_{ik} + \sum_{j \in I}
    \varepsilon_{ij}^{(v)} \log \frac{F_j^{v \to v_0}(T \cdot \mathbf{X}^{(v)}(g))}{F_j^{v \to v_0}(\mathbf{X}^{(v)}(g))}
\end{align*}
where $\mathbf{u}^{(v_0)}_g = (u^{(v_0)}_{g, i})_{i \in I}$.
By \eqref{eq:F_exp}, we have
\begin{align}
    \log \frac{F_j^{v \to v_0}(T \cdot \mathbf{X}^{(v)}(g))}{F_j^{v \to v_0}(\mathbf{X}^{(v)}(g))}
    &= \log \frac{1 + \sum_{\alpha} \xi_\alpha e^{\alpha_k} \cdot \prod_{\ell} X_\ell^{(v)}(g)^{\alpha_\ell} + e^{f_{jk}} \prod_{\ell} X_\ell^{(v)}(g)^{f_{j\ell}}}{1 + \sum_{\alpha} \xi_\alpha\cdot \prod_{\ell} X_\ell^{(v)}(g)^{\alpha_\ell} + \prod_{\ell} X_\ell^{(v)}(g)^{f_{j\ell}}}\label{eq:log_frac_F}.
\end{align}
On the other hand, we have
\begin{align*}
    X_\ell^{(v)}(g) = 
    \prod_{m \in I} X_m^{(v_0)}(g)^{\overline{c}_{\ell m}} \cdot F_m^{v_0 \to v}(\mathbf{X}^{(v_0)}(g))^{\varepsilon^{(v_0)}_{\ell m}}
\end{align*}
where $C^\bs_{v_0 \to v} = (\overline{c}_{ij})_{i,j \in I}$.
Moreover, by \cite[Proposition 5.3]{FZ-CA4}, we have
\begin{align*}
    F_{m}^{v_0 \to v}(\mathbf{X}^{(v_0)}(g)) = 
    \prod_{p \in I} X_p^{(v_0)}(g)^{\overline{f}_{mp}} \cdot \check{F}_m^{v_0 \to v}(\mathbf{X}^{(v_0)}(g)^{-1})
\end{align*}
where $F_{v_0 \to v}^\bs = (\overline{f}_{ij})_{i,j \in I}$ and 
$\check{F}^{v_0 \to v}_m$ be the $m$-th $F$-polynomial from $v_0$ to $v$ in the opposite mutation class $-\bs$.
By \cite[Theorem 2.8]{FuGy},
\begin{align*}
    C^\bs_{v_0 \to v} + \varepsilon^{(v_0)} F^\bs_{v_0 \to v} = C^{-\bs}_{v_0 \to v}.
\end{align*}
Therefore, 
\begin{align}\label{eq:X_asymp_infty}
    X_\ell^{(v)}(g) \sim
    \prod_{m \in I} X_m^{(v_0)}(g)^{\check{c}_{\ell m}} \quad \mbox{as } \mathbf{X}^{(v_0)}(g)^{-1} \to \mathbf{0}
\end{align}
where $C^{-\bs}_{v_0 \to v} = (\check{c}_{ij})_{i,j \in I}$.
Thus, we have
\begin{align*}
    &\log \frac{F_j^{v \to v_0}(T \cdot \mathbf{X}^{(v)}(g))}{F_j^{v \to v_0}(\mathbf{X}^{(v)}(g))}
    \sim \log \frac{1 + \sum_{\alpha} \xi_\alpha e^{\alpha_k} \cdot \prod_{m} X_m^{(v_0)}(g)^{\sum_\ell \alpha_\ell \check{c}_{\ell m}} + e^{f_{jk}} \prod_{m} X_m^{(v_0)}(g)^{\sum_\ell f_{j\ell} \check{c}_{\ell m}}}{1 + \sum_{\alpha} \xi_\alpha\cdot \prod_{m} X_m^{(v_0)}(g)^{\sum_\ell \alpha_\ell \check{c}_{\ell m}} + \prod_{m} X_m^{(v_0)}(g)^{\sum_\ell f_{j\ell} \check{c}_{\ell m}}}
\end{align*}
as $\mathbf{X}^{(v_0)}(g)^{-1} \to \mathbf{0}$.
Since there is a green path $\gamma: v_0 \to v$ for any $v \in \bExch_\bs$, the exponents $\sum_\ell \alpha_\ell \check{c}_{\ell m}$ and $\sum_\ell f_{j\ell} \check{c}_{\ell m}$ are negative by \cite[Corollary 5.5]{Gup}\footnote{In \cite{Gup}, the exponents $\sum_\ell \alpha_\ell \check{c}_{\ell m}$ and $\sum_\ell f_{j\ell} \check{c}_{\ell m}$ are the minus of the exponents of the \emph{generalized $F$-polynomial} \cite[Definition 2.25]{Gup} since $(C_{v \to v_0}^{\bs})^{-1} = C_{v_0 \to v}^{-\bs}$. The positivity of the exponents of the generalized $F$-polynomial is equivalent with it is a polynomial in fact. She state that the existence of the green sequence is a sufficient condition of the positivity in \cite[Corollary 5.5]{Gup}.} and the existence of the green paths $\gamma: v \to v_0$ in $\bExch_\bs$ (\cite{AIR}, see also \cite[Theorem 3.6]{BY}).
Hence,
\begin{align*}
    \lim_{\mathbf{X}^{(v_0)}(g)^{-1} \to \mathbf{0}} \log \frac{F_j^{v \to v_0}(T \cdot \mathbf{X}^{(v)}(g))}{F_j^{v \to v_0}(\mathbf{X}^{(v)}(g))} = 0.
\end{align*}
Therefore, 
\begin{align*}
    \lim_{\mathbf{X}^{(v_0)}(g)^{-1} \to \mathbf{0}} 
    \big(u^{(v_0)}_{g,i}(L^{(v_0)}_k) \big)_{i,k \in I}
    = C_{v \to v_0}^\bs
    = (G_{v_0 \to v}^{-\bs})^\tr.
\end{align*}
By \cref{prop:cl_cpx=g-vect_fan}, we get $\lim_{\mathbf{X}^{(v_0)}(g)^{-1} \to \mathbf{0}} \mathbf{u}^{(v_0)}_{g}(\fF^+_{\bs}) = \mathbf{x}^{(v_0)}(\fF^+_{\bs})$.
\end{proof}

This theorem only describes the asymptotic behavior of the simplicial set $\mathbf{u}^{(v_0)}_{g}(\fF^+_{\bs})$ in the direction  $\interior \bS \cC^\pm_{(v_0)}$.
\begin{prob}
Describe the asymptotic behavior of $\mathbf{u}^{(v_0)}_{g}(\fF^+_{\bs})$ when $g$ goes toward the boundary $\bS \cX_\bs(\bR^\trop)$ except the subset $\interior \bS \cC^+_{(v_0)}$.
\end{prob}

\begin{proof}[Proof of \cref{thm:fan_in_tangent}]
From the definition, we have
\begin{align*}
    (d \mathbf{u}^{(v_0)}_{g})_0 = (d \log \mathbf{X}^{(v_0)})_{g} \circ (d E(g,-))_0
\end{align*}
for $g \in \cX_\bs(\pos)$ and $v_0 \in \bExch_\bs$.
Thus, this statement directly follows from \cref{thm:g_to_infty}.
\end{proof}

\appendix

\section{Cluster algebras of surface type}\label{cl_surf}
In this subsection, we review mutation classes of surface type. For basic notions in hyperbolic geometry, we refer the reader to \cite{CB,Penner}. 

Let $\Sigma$ be an oriented compact topological surface which may have boundaries and $M$ be a finite set in $\Sigma$. We call the points of $M$ \emph{marked points}, and  $(\Sigma,M)$ a \emph{marked surface}. 
Let $\interior \Sigma$ and $\partial \Sigma$ denote the interior and the boundary of $\Sigma$, respectively. We call a point of $M_\circ:=M\cap \interior\Sigma$ a \emph{puncture}, and a point of $M_\partial:=M \cap \partial\Sigma$ a \emph{special point}. Let $\Sigma^\ast:=\Sigma \setminus M_\circ$. 
We assume the two following conditions:
\begin{itemize}
    \item each boundary component has at least one special point;
    \item $-2\chi(\Sigma^\ast) + |M_\partial|>0$.\footnote{The left-hand side counts the number of triangles of an ideal triangulation of $\Sigma$. In particular, this condition excludes the bigon, which does not admit any ideal triangulation.}
\end{itemize}
An \emph{ideal arc} in $(\Sigma,M)$ is an isotopy class of arc in $\Sigma^\ast$ connecting marked points without self-intersection except for its endpoints. In this paper, we require an ideal arc not to be represented by a boundary segment connecting two consecutive special points. 
A maximal disjoint collection $\tri$ of ideal arcs is called an \emph{ideal triangulation} of $\Sigma$.
Let $n \ceq -3 \chi(\Sigma^*) + |M_\partial|$ be the number of ideal arcs. A \emph{labeled triangulation} is an ideal triangulation $\tri$ equipped with a bijection $\ell: I \to \tri$, where $I=\{1,\dots,n\}$ is a fixed index set. 

\begin{figure}[ht]
    \centering
\begin{tikzpicture}
\draw[blue](0,0) -- (0,1.8);
\draw[blue] (0,0) ..controls (45:0.5) and (1.5,2.5).. (0,2.5);
\draw[blue] (0,0) ..controls (135:0.5) and (-1.5,2.5).. (0,2.5);
\fill(0,0) circle(2pt);
\fill(0,1.8) circle(2pt);
\end{tikzpicture}
    \caption{A self-folded triangle.}
    \label{fig:self-folded}
\end{figure}

\begin{defi}
    For a labeled triangulation $(\tri,\ell)$ of $\Sigma$ without \emph{self-folded triangles} as shown in \cref{fig:self-folded}, 
we define a \emph{quiver} (\emph{i.e.}, a directed graph) $Q^{(\tri,\ell)}$ with the vertex set $I$ as follows. 
If a triple $(\ell(i),\ell(j),\ell(k))$ compose an ideal triangle of $\tri$, then we connect any pair of them by a counterclockwise arrow on that triangle. A local picture is shown in \cref{fig:quiv_tri}. 
Then let $\ve^{(\tri,\ell)}=(\ve^{(\tri,\ell)}_{ij})_{i,j \in I}$ be the corresponding $n\times n$ skew-symmetric exchange matrix, which is given by 
\[
    \ve_{ij}^{(\tri,\ell)}:= \#\{\text{arrows from $i$ to $j$}\} - \#\{\text{arrows from $j$ to $i$}\}.
\]
For the case with self-folded triangles, we assign the matrix according to \cite[Definition 4.1]{FST08}.
\end{defi}

\begin{figure}[h]
    \begin{align*}
    \tikz[scale=1.2,>=latex]{
    \path(0,0) node [fill, circle, inner sep=1.5pt] (x1){};
    \path(135:2) node [fill, circle, inner sep=1.5pt] (x2){};
    \path(0,2*1.4142) node [fill, circle, inner sep=1.5pt] (x3){};
    \path(45:2) node [fill, circle, inner sep=1.5pt] (x4){};
    \draw[blue](x1) to (x2) to (x3) to (x4) to (x1) to  (x3);
    \node[blue] at (1,2.75) {$\tri$};
    \node[mygreen] at (1,0) {$Q^\tri$};
    \color{mygreen}{
        \draw(0,1.4142) circle(2pt) coordinate(v0);
        \draw(135:2)++(45:1) circle(2pt) coordinate(v1);
        \draw(45:2)++(135:1) circle(2pt) coordinate(v2);
        \draw(45:1) circle(2pt) coordinate(v3);
        \draw(135:1) circle(2pt) coordinate(v4);
        \qarrow{v0}{v4}
        \qarrow{v4}{v1}
        \qarrow{v1}{v0}
        \qarrow{v0}{v2}
        \qarrow{v2}{v3}
        \qarrow{v3}{v0}
        }
    }
    \end{align*}
    \caption{The quiver $Q^\tri$ associated with an ideal triangulation $\tri$ of $\Sigma$.}
    \label{fig:quiv_tri}
\end{figure}
It is classically known that any two ideal triangulations $\tri,\tri'$ of a common marked surface $\Sigma$ can be connected by a finite sequence of flips, as shown in \cref{fig:flip} (see, for example, \cite{Pen87}). Using the ``common'' labeling for them, one can verify that the associated exchange matrices $\ve^{(\tri,\ell)}$ and $\ve^{(\tri',\ell)}$ are related by the corresponding sequence of matrix mutations. The exchange matrices $\ve^{(\tri,\ell)}$, $\ve^{(\tri,\ell')}$ with different labelings are obviously related by a sequence of transpositions. 
In this way, we get a canonical mutation class $\bs_\Sigma$ of exchange matrices $\ve^{(\tri,\ell)}$ which only depends on $\Sigma$. Let us denote the corresponding objects by
\begin{align*}
    \bExch_\Sigma:=\bExch_{\bs_\Sigma},\quad \X_\Sigma(\pos):=\X_{\bs_\Sigma}(\pos), \quad \mbox{and}\quad \X_\Sigma(\bR^\trop):=\X_{\bs_\Sigma}(\bR^\trop).
\end{align*}
Consider the graph $\mathbb{T}\mathrm{ri}_\Sigma$ with the vertex set given by labeled triangulations of $\Sigma$, and with the labeled edges of the following types:
\begin{itemize}
    \item a labeled edge of the form $(\tri,\ell) \overbar{k} (\tri',\ell)$ corresponding to the flip along $\ell(k) \in \tri$;
    \item a labeled edge edge of the form $(\tri,\ell) \overbarnear{\sigma} (\tri,\ell')$ corresponding to the transposition $\sigma=(j\ k)$ of labelings for $(j,k) \in I \times I$.
\end{itemize}
Then by construction, we have a full graph embedding
\begin{align}\label{eq:graph_emb}
  \mathbb{T}\mathrm{ri}_\Sigma \hookrightarrow \bExch_\Sigma.
\end{align}
In general, we can use the \emph{tagged triangulations} \cite{FST08} to realize all the vertices in $\bExch_\Sigma$.
Through this graph isomorphism, we can get the isomorphism between the extended mapping class group $MC(\Sigma) \ltimes (\bZ/2)^{M_\circ}$ and the cluster modular group $\Gamma_{\bs_\Sigma}$ of $\bs_\Sigma$ except for few cases (cf. \cite{BS}).

\begin{figure}[h]
    \centering
    \begin{tikzpicture}[scale=0.7, >=latex]
    \path(0,0) node [fill, circle, inner sep=1.6pt] (x1){};
    \path(135:4) node [fill, circle, inner sep=1.6pt] (x2){};
    \path(0,4*1.4142) node [fill, circle, inner sep=1.6pt] (x3){};
    \path(45:4) node [fill, circle, inner sep=1.6pt] (x4){};
    \draw[blue](x1) to node[midway,left,black]{} (x2) 
    to node[midway,left,black]{} (x3) 
    to node[midway,right,black]{} (x4) 
    to node[midway,right,black]{} (x1) 
    to node[midway,left,black]{} (x3);
        {\color{mygreen}
        \draw(0,2*1.4142) circle(3pt) coordinate(v0);
        \draw(135:4)++(45:2) circle(3pt) coordinate(v1);
        \draw(45:4)++(135:2) circle(3pt) coordinate(v2);
        \draw(45:2) circle(3pt) coordinate(v3);
        \draw(135:2) circle(3pt) coordinate(v4);
        \qarrow{v0}{v4}
        \qarrow{v4}{v1}
        \qarrow{v1}{v0}
        \qarrow{v0}{v2}
        \qarrow{v2}{v3}
        \qarrow{v3}{v0}
        }
    
    \draw[-implies, double distance=2pt](4,2*1.4142) to node[midway,above]{$f_{k}$} (6,2*1.4142);
    
    \begin{scope}[xshift=10cm]
    \path(0,0) node [fill, circle, inner sep=1.6pt] (x1){};
    \path(135:4) node [fill, circle, inner sep=1.6pt] (x2){};
    \path(0,4*1.4142) node [fill, circle, inner sep=1.6pt] (x3){};
    \path(45:4) node [fill, circle, inner sep=1.6pt] (x4){};
    \draw[blue](x1) to node[midway,left,black]{} (x2) 
    to node[midway,left,black]{} (x3) 
    to node[midway,right,black]{} (x4) 
    to node[midway,right,black]{} (x1);
    \draw[blue] (x2) to node[midway,above,black]{} (x4);
        {\color{mygreen}
        \draw(0,2*1.4142) circle(3pt) coordinate(v0);
        \draw(135:4)++(45:2) circle(3pt) coordinate(v1);
        \draw(45:4)++(135:2) circle(3pt) coordinate(v2);
        \draw(45:2) circle(3pt) coordinate(v3);
        \draw(135:2) circle(3pt) coordinate(v4);
        \qarrow{v0}{v1}
        \qarrow{v1}{v2}
        \qarrow{v2}{v0}
        \qarrow{v0}{v3}
        \qarrow{v3}{v4}
        \qarrow{v4}{v0}
        }
    \end{scope}
    \node [blue] at (0.25,3.5) {$k$};
    \node [blue] at (10.75,2.5) {$k$};

    \end{tikzpicture}
    \caption{Flip along the edge $\ell(k)$. Here we give the same label $k$ to the resulting edge.}
    \label{fig:flip}
\end{figure}

In the following, we review the geometric models for the spaces $\X_\Sigma(\pos)$ and $\X_\Sigma(\bR^\trop)$.

\smallskip
\paragraph{\textbf{The cluster $\X$-manifold $\X_\Sigma(\pos)$ and the enhanced \Teich\ space.}}
Let $\HHH$ denote the Poincar\'e upper-half plane model of the hyperbolic plane, and $\SSS=\bR \cup \{\infty\}$ its ideal boundary. 

\begin{defi}
     Let $\mathrm{Hyp}(\Sigma)$ be the set of hyperbolic metrics on $\Sigma \setminus M$ where each point in 
     each puncture
     corresponds to a funnel or a cusp, and 
     each special point corresponds to a spike of a crown. For $h \in \mathrm{Hyp}(\Sigma)$, let $M_\circ^h \subset M_\circ$ denote the set of punctures that give rise to funnels. Let $\mathrm{Hyp}_0(\Sigma) \subset \mathrm{Hyp}(\Sigma)$ be the subset consisting of those without cusps. 
\end{defi}
See \cref{fig:Sigma^h} for an example.


\begin{defi}
    \begin{enumerate}\label{def_Teich}
        \item For $h\in{\rm Hyp}(\Sigma)$, let $\Sigma^h$ be the hyperbolic surface obtained from $(\Sigma \setminus M,h)$ by truncating the outer side of the shortest closed geodesic in each funnel. See \cref{fig:Sigma^h}. For $p \in M_\circ^h$, let $\partial_p$ denote the resulting boundary component. 
        \item \label{def_sig} We call an orientation-preserving homeomorphism $f\cl \interior (\Sigma \setminus M)\rightarrow \interior \Sigma^h$ a \emph{signed homeomorphism} if it maps a representative of each ideal arc to a complete geodesic so that each end incident to $p \in M_\circ^h$ is mapped to a spiralling geodesic into the geodesic boundary $\partial_p$ in either left or right direction. 
        Given a signed homeomorphism $f$, its \emph{signature} $\epsilon_f=(\eta_p)_p\in\{+,0,-\}^{M_\circ}$ is defined by 
        \[
            \eta_p=\begin{cases}
            + &\mbox{if $p\in M_\circ^h$ and $f(\alpha_p)$ spirals to the left along the geodesic boundary $\partial_p$},\\
                0&\mbox{if $p \in M_\circ \setminus M_\circ^h$}, \\
                - &\mbox{if $p\in M_\circ^h$ and $f(\alpha_p)$ spirals to the right along the geodesic boundary $\partial_p$}.
            \end{cases}
        \]
        Here $\alpha_p$ is any ideal arc incident to $p$. 
        \item We say that two pairs $(\Sigma^{h_i},f_i)$ for $i=1,2$ are \emph{equivalent} if $f_2\circ f_1^{-1}$ is homotopic to an isometry $\Sigma^{h_1} \to \Sigma^{h_2}$, and $\epsilon_{f_1}=\epsilon_{f_2}$. We call the set $\widehat{\cT}(\Sigma)$ of equivalence classes of the pairs $(\Sigma^h,f)$ the \emph{enhanced \Teich\ space}. 
    \end{enumerate}
\end{defi}

\begin{figure}[h]
    \centering
    \begin{tikzpicture}[scale=.7]
    \draw  (-4.5,0) ellipse (3 and 2);
    \draw (-5.5,-0.5) .. controls (-5.5,-1.35) and (-3.5,-1.35) .. (-3.5,-0.5);
    \draw (-5.4,-0.8) .. controls (-5.4,-0.2) and (-3.6,-0.2) .. (-3.6,-0.8);
    \draw  (-5.5,0.65) ellipse (0.5 and 0.5);
    \node [fill, circle, inner sep=1.3pt] at (-3.5,0.65) {};
    \node [fill, circle, inner sep=1.3pt] at (-5.05,0.85) {};
    \node [fill, circle, inner sep=1.3pt] at (-5.95,0.85) {};
    \node [fill, circle, inner sep=1.3pt] at (-5.5,0.15) {};
    \node at (-3.1,0.6) {$p$};
    \draw[red] (-5.5,0.15) .. controls (-5.2,-0.35) and (-3.5,0) .. (-3.5,0.65);
    \node [red] at (-4.4,0.25) {$\alpha$};

    \draw  (3,0) ellipse (3 and 2);
    \draw [white, ultra thick](1.3,1.65) .. controls (1.6,1.8) and (2,1.9) .. (2.3,1.95);
    \draw [white, ultra thick](3.75,1.95) .. controls (4,1.9) and (4.3,1.8) .. (4.55,1.7);
    \draw (2,-0.5) .. controls (2,-1.35) and (4,-1.35) .. (4,-0.5);
    \draw (2.1,-0.8) .. controls (2.1,-0.2) and (3.9,-0.2) .. (3.9,-0.8);
    \draw (0.85,0.7) .. controls (1.35,0.95) and (1.35,2.3) .. (1.35,3.15) .. controls (1.35,2.3) and (1.8,1.8) .. (1.8,2.65) .. controls (1.8,1.8) and (2.25,2.3) .. (2.25,3.15) .. controls (2.25,2.3) and (1.35,2.3) .. (1.35,3.15);
    \draw (2.75,0.7) .. controls (2.25,0.95) and (2.25,2.3) .. (2.25,3.15);
    \draw (3.2,0.7) .. controls (3.7,0.95) and (3.75,1.95) .. (3.75,2.5);
    \draw (5.05,0.7) .. controls (4.55,0.95) and (4.5,1.95) .. (4.5,2.5);
    \draw [dashed](3.75,2.5) .. controls (3.75,2.75) and (3.75,2.9) .. (3.6,3.2);
    \draw [dashed](4.5,2.5) .. controls (4.5,2.75) and (4.5,2.9) .. (4.65,3.2);

    \draw  (10.5,0) ellipse (3 and 2);
    \draw [white, ultra thick](8.8,1.65) .. controls (9.1,1.8) and (9.5,1.9) .. (9.8,1.95);
    \draw [white, ultra thick](11.25,1.95) .. controls (11.5,1.9) and (11.8,1.8) .. (12.05,1.7);
    \draw (9.5,-0.5) .. controls (9.5,-1.35) and (11.5,-1.35) .. (11.5,-0.5);
    \draw (9.6,-0.8) .. controls (9.6,-0.2) and (11.4,-0.2) .. (11.4,-0.8);
    \draw (8.35,0.7) .. controls (8.85,0.95) and (8.85,2.3) .. (8.85,3.15) .. controls (8.85,2.3) and (9.3,1.8) .. (9.3,2.65) .. controls (9.3,1.8) and (9.75,2.3) .. (9.75,3.15) .. controls (9.75,2.3) and (8.85,2.3) .. (8.85,3.15);
    \draw (10.25,0.7) .. controls (9.75,0.95) and (9.75,2.3) .. (9.75,3.15);
    \draw (10.7,0.7) .. controls (11.2,0.95) and (11.25,1.95) .. (11.25,2.5);
    \draw (12.55,0.7) .. controls (12.05,0.95) and (12.05,1.95) .. (12.05,2.5);
    \draw(11.65,2.5) ellipse (0.4 and 0.2);
    \node at (12.5,2.5) {$\partial_p$};
    
    \node at (-4.5,-2.5) {$\Sigma$};
    \node at (3,-2.5) {$(\Sigma \setminus M, h)$};
    \node at (10.5,-2.5) {$\Sigma^h$};
    \node [red] at (12,.3) {$f(\alpha)$};
    \draw [red](9.3,2.65) .. controls (9.3,0.5) and (10,0.25) .. (10.5,0.25) .. controls (11,0.25) and (11.65,0.5) .. (11.65,1.35) .. controls (11.65,1.7) and (11.45,1.9) .. (11.23,2);
    \draw [red](12.05,2.05) .. controls (11.85,1.95) and (11.45,1.95) .. (11.25,2.15);\draw [red](12.05,2.25) .. controls (11.85,2.1) and (11.45,2.1) .. (11.25,2.3);
    \draw [red](12.05,2.4) .. controls (11.85,2.2) and (11.45,2.2) .. (11.25,2.4);
    \end{tikzpicture}
    \caption{Left: a topological surface $\Sigma$ and an ideal arc $\alpha$ (red arc); center: the surface $\Sigma \setminus M$ with a hyperbolic metric $h$ such that $p \in M^h_\circ$; right: the surface $\Sigma^h$ and the arc $f(\alpha)$ if $\eta_p = +$.}
    \label{fig:Sigma^h}
\end{figure}

\begin{defi}\label{cross_ratio}
    Given a labeled triangulation $(\tri,\ell)$ and $i \in I$, let $\square_i$ be the unique quadrilateral in $\tri$ containing the ideal arc $\ell(i)$ as a diagonal. For a point $[\Sigma^h,f] \in \widehat{\cT}(\Sigma)$, choose a geodesic lift $\widetilde{f(\square_i)} \subset \widetilde{\Sigma^h} \subset \HHH$ of $f(\square_i)$ with respect to the universal covering $\widetilde{\Sigma^h} \to \Sigma^h$ determined by the hyperbolic structure $h$. Let $x,y,z,w \in \SSS$ be the four ideal vertices of $\widetilde{f(\square_i)}$ in this counter-clockwise order, and assume that $x$ is an endpoint of the lift of $\ell(i)$. 
    Then we define the \emph{cross ratio coordinate}  associated with $i \in I$ by
    \[X_i^{(\tri,\ell)}([\Sigma^h,f])
    \ceq [x : y : z : w] >0,
    \]
    where 
    \begin{align*}
        [x:y:z:w]&\ceq \frac{w-x}{w-z}\cdot\frac{z-y}{y-x} 
    \end{align*}
    denotes the cross ratio of four points in $\SSS$.  
    Since the right-hand side is invariant under M\"{o}bius transformations, $X_i^{(\tri,\ell)}$ does not depend on the choice of the lift $\widetilde{f(\square_i)}$, and does not change if we choose the other endpoint of $\ell(i)$ as $x$.
\end{defi}

\begin{thm}[Fock--Thurston (cf.~\cite{FG07})]\label{shear}
    For any labeled triangulation $(\tri,\ell)$ of $\Sigma$,
    the map
    \[
        \mathbf{X}^{(\tri,\ell)}: \widehat{\cT}(\Sigma)\rightarrow\bR^{I}_{>0},\quad [\Sigma^h,f]\mapsto (X_i^{(\tri,\ell)}([\Sigma^h,f]))_{i\in I}
    \]
    is bijective.
\end{thm}

It is also known that the coordinate transformation $\mathbf{X}^{(\tri',\ell)} \circ (\mathbf{X}^{(\tri,\ell)})^{-1}$ associated with a flip coincides with the cluster $\X$-transformation \cite{FG07}. 
Therefore, the pairs $(\ve^{(\tri,\ell)},\mathbf{X}^{(\tri,\ell)})$ associated with labeled triangulations $(\tri,\ell)$ generates a canonical mutation class $\bs_\Sigma$ in the ambient field $\cF_X:=\Frac\cC^\infty(\widehat{\cT}(\Sigma))$. Then \cref{shear} can be rephrased that we have a canonical isomorphism $\widehat{\cT}(\Sigma) \cong \X_{\Sigma}(\pos)$ that is compatible with \eqref{eq:graph_emb}. 
We say that a mutation class is \emph{of surface type} if it is isomorphic to the mutation class $\bs_\Sigma$ for some marked surface $\Sigma$.

\smallskip
Here is another description of $\cT(\Sigma)$, which we will use to define the earthquakes.
Let $\cF_\infty(\Sigma)$ be the \emph{Farey set} of the marked surface $\Sigma$, which is a cyclically ordered set defined as follows \cite[Section 1.3]{FG03}. Let $h_0 \in \mathrm{Hyp}(\Sigma)$ be an auxiliary complete hyperbolic metric on $\Sigma \setminus M$ of finite area. Then the universal cover $\widetilde{\Sigma} \subset \HHH$ is a convex set, and $\cF_\infty(\Sigma) \subset \SSS$ is defined to be the set of its ideal vertices (\emph{i.e.}, lifts of points in $M$). It inherits the cyclic ordering of $\SSS$. Note that, as a cyclically ordered set, the Farey set $\cF_\infty(\Sigma)$ does not depend on $h_0$. 
\begin{prop}[cf.~{\cite[Theorem 1.6]{FG03}}]\label{prop:enhanced_Farey}
    There is a natural bijection between the enhanced \Teich\ space $\widehat{\cT}(\Sigma)$ and the space of $\PSL$-orbits of pairs $(\rho,\psi)$, where $\rho:\pi_1(\Sigma) \to \PSL$ is a Fuchsian representation and $\psi: \cF_\infty(\Sigma) \to \SSS$ is a $\rho$-equivariant, order-preserving map. 
\end{prop}
Indeed, given $[\Sigma^h,f] \in \widehat{\cT}(\Sigma)$, the representation $\rho$ is the monodromy of the hyperbolic structure $h$, and $\psi$ is given by the continuous extension of the lift $\widetilde{f}:\widetilde{\Sigma} \to \widetilde{\Sigma^h} \subset \HHH$ to the ideal boundary. 
\smallskip

\paragraph{\textbf{The tropical $\X$-manifold $\X_\Sigma(\bR^\trop)$ and the space of measured geodesic laminations.}}
For basics on the measured geodesic laminations, we refer the reader to \cite{PH,BKS16}. 

We fix a reference metric $h\in\mathrm{Hyp}(\Sigma)$. 
Recall that a \emph{geodesic lamination} is a closed subset $G$ in $\Sigma^h$ which is a disjoint union of complete simple geodesics. Each geodesic in a geodesic lamination is called a \emph{leaf}. 
The following lemma describes the behavior of a geodesic lamination near each closed geodesic boundary. 

\begin{lem}[{\cite[Lemma 2.1]{BKS16}}]\label{lem:spiral}
    For $p \in M_\circ^h$, there exists a neigborhood $U$ of  $\partial_p$ such that each leaf in $G\cap U$ spirals along $\partial_p$ either in the left or right direction (i.e., the direction following or against the orientation of $\partial_p$ induced by $\Sigma$).
\end{lem}
Let us define the \emph{lamination signature} $\sigma_G =(\sigma_p)_p \in \{+,0,-\}^{M_\circ^h}$ by $\sigma_p:=+$ (resp. $\sigma_p:=-$) if $G$ spirals along $\partial_p$ in the left (resp. right) direction and by $\sigma_p:=0$ otherwise.

\begin{lem}[{\cite[Lemma 2.3]{BKS16}}]\label{lem:spiral_modify}
For an ideal arc $e$ connecting punctures $p_1,p_2 \in M_\circ^h$ and any signs $\sigma_1,\sigma_2 \in \{+,-\}$, there exists a unique spiralling geodesic $G$ on $\Sigma^h$ homotopic to $e$ such that $\sigma_{p_i}(G)=\sigma_i$ for $i=1,2$. 
\end{lem}
A \emph{tranverse measure} on a geodesic lamination $G$
is an assignment
of each \emph{transverse arc} $\alpha$ (\emph{i.e.}, an arc which intersects $G$ transversely and whose endpoints are on $\Sigma^h\setminus G$)
to $\mu(\alpha) \in\bR_{>0}$
such that $\sigma$-additive for concatenation of arcs, invariant under the isotopy relative to the leaves of $G$
, and that $\mu(\alpha)=0$ if and only if $\alpha \cap G = \emptyset$.  
We call a pair $(G,\mu)$ a \emph{measured geodesic lamination}. The set of the measured geodesic laminations on $\Sigma^h$ is denoted by $\ML(\Sigma^h)$. When we study the earthquake maps, it is useful to combine these spaces to build a bundle over the enhance \Teich\ space $\widehat{\cT}(\Sigma)$. 
    
\begin{defi}[cf.~\cite{BB09}]
Let $[\Sigma^h,f] \in \widehat{\cT}(\Sigma)$. 
For $(G,\mu) \in \ML(\Sigma^h)$, a \emph{relaxed signature under the marking $f$} is a tuple $\eta=(\eta_p)_p \in\{+,0,-\}^{M_\circ}$ such that $\eta|_{M_\circ^h} = \sigma_G\cdot \epsilon_f$, and $\eta_p =0$ (resp.~$\eta_p \in \{+,-\}$) if $p \in M_\circ \setminus M_\circ^h$ and there is no leaf (resp.~a leaf) incident to the corresponding cusp. 
    Then we define 
\[
    \widehat{\ML}([\Sigma^h,f])\ceq\{(G,\mu,\eta)\mid (G,\mu)\in\ML(\Sigma^h)\text{, $\eta$: a relaxed signature of $G$ under $f$}\}.
\]
Then the \emph{bundle of measured geodesic lamination} is a fiber bundle $\widehat{\ML}_\Sigma \to \widehat{\cT}(\Sigma)$ such that the fiber over a point $[\Sigma^h,f] \in \widehat{\cT}(\Sigma)$ is given by $\widehat{\ML}([\Sigma^h,f])$. 
\end{defi}

We are going to define tropical $\cX$-variables on each fiber $\eML([\Sigma^h,f])$ for $[\Sigma^h,f] \in \widehat{\cT}(\Sigma)$. We call an orientation-preserving homeomorphism $f^\bot: \interior\Sigma \to \interior\Sigma^h$ a \emph{perpendicular homeomorphism} if it maps a representative of each ideal arc to a geodesic so that each end incident to $p \in M_\circ^h$ is mapped to a perpendicular geodesic to the geodesic boundary $\partial_p$. Notice that for each signed homeomorphism $f: \interior\Sigma \to \interior\Sigma^h$, there exists a unique perpendicular homeomorphism $f^\bot: \interior\Sigma \to \interior\Sigma^h$ up to homotopy such that the composition $f^\bot \circ f^{-1}|_{\interior\Sigma^h \setminus U}$ preserves the based homotopy class of each properly embedded arc. Here $U$ is the union of neighborhoods of $\partial_p$ for all $p \in M_\circ^h$ as in \cref{lem:spiral}.

First we consider the case $h \in \mathrm{Hyp}_0(\Sigma)$, until just before \cref{shear1}.
Fix a labeled triangulation $(\tri,\ell)$ of $\Sigma$, and represent each ideal arc by perpendicular geodesics $f^\bot(\ell(i))$ for $i \in I$. Given $(G,\mu,\eta) \in \widehat{\ML}([\Sigma^h,f])$, let $(G',\mu')$ be the measured geodesic lamination obtained from $(G,\mu)$ by:
\begin{itemize}
    \item changing the spiralling direction of leaves of $G$ so that $\sigma_{G'}=\eta=\sigma_G \cdot \epsilon_f$ (It is possible by \cref{lem:spiral_modify}), and 
    \item shifting each leaf of $G'$ incident to a spike to a geodesic transverse to the boundary following the orientation induced by $\Sigma$, as follows:
\end{itemize}

\begin{center}
\begin{tikzpicture}[scale=.8]
\draw (-1.5,2) .. controls (-1.5,0.55) and (0.5,0.55) .. (0.5,2);
\draw [white, ultra thick](-0.6,0.9)  -- (-0.4,0.9);
\draw (-3,-0.5) .. controls (-2,0) and (-1.5,0.55) .. (-1.5,2) .. controls (-1.5,0.55) and (-0.5,0.05) .. (-0.5,1.5) .. controls (-0.5,0.05) and (0.5,0.55) .. (0.5,2) .. controls (0.5,0.55) and (1,0) .. (2,-0.5);
\draw [red](-0.5,1.5) .. controls (-0.5,0.05) and (0,-0.1) .. (0.25,-0.25);
\draw [red, dashed](0.25,-0.25) .. controls (0.5,-0.4) and (0.75,-0.45) .. (1,-0.5);

\draw [very thick,-{Classical TikZ Rightarrow[length=4pt]},decorate,decoration={snake,amplitude=1.5pt,pre length=2pt,post length=3pt}](2.5,0.75) -- (4,0.75);

\draw (6,2) .. controls (6,0.55) and (8,0.55) .. (8,2);
\draw [white, ultra thick](6.9,0.9)  -- (7.1,0.9);
\draw (4.5,-0.5) .. controls (5.5,0) and (6,0.55) .. (6,2) .. controls (6,0.55) and (7,0.05) .. (7,1.5) .. controls (7,0.05) and (8,0.55) .. (8,2) .. controls (8,0.55) and (8.5,0) .. (9.5,-0.5);
\draw [red](6.55,0.65) .. controls (7,0.05) and (7.5,-0.15) .. (7.75,-0.25);
\draw [red, dashed](7.75,-0.25) .. controls (8,-0.35) and (8.25,-0.45) .. (8.5,-0.5);
\end{tikzpicture}
\end{center}
Then, the resulting ``shifted'' measured geodesic lamination $(G',\mu')$ transversely intersects with geodesics $f^\bot(\ell(i))$.

For each $i \in I$, let $\square_i$ be the unique quadrilateral in $\tri$ which contains $i$ as a diagonal. Let $p_1,p_2 \in M_\circ^h$ be the punctures connected by the ideal arc $\ell(i)$. 
Fix a parametrization $\bR \xrightarrow{\sim} f^\bot(\ell(i))$ of the perpendicular geodesic. 
Let us draw the picture so that $f^\bot(\ell(i))$ is a vertical edge and $\partial_{p_1}$ lies on the top: see \cref{fig:lam_in_quad}. 
Then each path-connected component $\gamma$ of $G'\cap \square_i$ satisfies either one of the followings: 
(i) $\gamma$ intersects the upper-left and lower-right sides of $\square_i$;
(ii) $\gamma$ intersects the lower-left and upper-right sides of $\square_i$;
(iii) othewise. 
Here path-connected components satisfying (i) and (ii) cannot coexist in $\square_i$. 

Consider the closed subsets
\begin{align*}
    &S_{\mathrm{in}}:=\{ f^\bot(\ell(i)) \cap \gamma\mid \text{$\gamma$ is a component of type (i) or (ii)}\}, \\
    &S_{p_i}:=\{f^\bot(\ell(i)) \cap \gamma \mid \text{$\gamma$ is a component of type (iii), surrounding $\partial_{p_i}$}\}
\end{align*}
of $f^\bot(\ell(i))$, and consider $i_+:=(\min S_{p_1}+ \max S_{\mathrm{in}})/2$ and $i_-:=(\max S_{\mathrm{out}}+ \min S_{p_2})/2$. 
Then $\alpha_i:=[i_-,i_+] \subset f(\ell(i))$ is a transverse arc to $(G',\mu')$. Define
\[
    x_i^{(\tri,\ell)}(G,\mu) \ceq \begin{cases}
        \mu'(\alpha_i) & \mbox{if a component $\gamma$ of type (i) exists}, \\
        -\mu'(\alpha_i)& \mbox{if a component $\gamma$ of type (ii) exists}, \\
        0 & \mbox{otherwise}. 
    \end{cases}
\]

\begin{figure}[h]
    \centering
    \begin{tikzpicture}[scale=.8]
    \draw [blue](-1.5,3) coordinate (v1) -- (-4,0.5) node (v2) {} -- (-1.5,-2) node (v3) {} -- (1,0.5) node (v4) {} -- (v1);
    \draw[blue] (v1) -- (v3);
    \draw[fill=white]  (v1) ellipse (0.5 and 0.5);
    \draw[fill=white]  (v2) ellipse (0.5 and 0.5);
    \draw[fill=white]  (v3) ellipse (0.5 and 0.5);
    \draw[fill=white]  (v4) ellipse (0.5 and 0.5);
    \draw [red](-2.5,3) .. controls (-2.5,2.5) and (-2,2) .. (-1.5,2) .. controls (-1,2) and (-0.5,2.5) .. (-0.5,3);
    \draw [red](-4,1.5) .. controls (-3.5,1.5) and (-3,1) .. (-3,0.5) .. controls (-3,0) and (-3.5,-0.5) .. (-4,-0.5);
    \draw [red](-2.5,-2) .. controls (-2.5,-1.5) and (-2,-1) .. (-1.5,-1) .. controls (-1,-1) and (-0.5,-1.5) .. (-0.5,-2);
    \draw [red](1,-0.5) .. controls (0.5,-0.5) and (0,0) .. (0,0.5) .. controls (0,1) and (0.5,1.5) .. (1,1.5);
    \draw [red](-2.05,3) .. controls (-2.05,2.7) and (-1.8,2.45) .. (-1.5,2.45) .. controls (-1.2,2.45) and (-0.95,2.7) .. (-0.95,3);
    \draw [red](-4,1.05) .. controls (-3.7,1.05) and (-3.45,0.8) .. (-3.45,0.5) .. controls (-3.45,0.2) and (-3.7,-0.05) .. (-4,-0.05);
    \draw [red](-2.05,-2) .. controls (-2.05,-1.7) and (-1.8,-1.45) .. (-1.5,-1.45) .. controls (-1.2,-1.45) and (-0.95,-1.7) .. (-0.95,-2);
    \draw [red](1,-0.05) .. controls (0.7,-0.05) and (0.45,0.2) .. (0.45,0.5) .. controls (0.45,0.8) and (0.7,1.05) .. (1,1.05);
    \draw [red](-2.15,3) .. controls (-2.15,2.65) and (-1.85,2.35) .. (-1.5,2.35) .. controls (-1.15,2.35) and (-0.85,2.65) .. (-0.85,3);
    \draw [red](-4,1.15) .. controls (-3.65,1.15) and (-3.35,0.85) .. (-3.35,0.5) .. controls (-3.35,0.15) and (-3.65,-0.15) .. (-4,-0.15);
    \draw [red](-2.15,-2) .. controls (-2.15,-1.65) and (-1.85,-1.35) .. (-1.5,-1.35) .. controls (-1.15,-1.35) and (-0.85,-1.65) .. (-0.85,-2);
    \draw [red](1,-0.15) .. controls (0.65,-0.15) and (0.35,0.15) .. (0.35,0.5) .. controls (0.35,0.85) and (0.65,1.15) .. (1,1.15);
    \draw [red](-2.3,3) .. controls (-2.3,2.55) and (-1.95,2.2) .. (-1.5,2.2) .. controls (-1.05,2.2) and (-0.7,2.55) .. (-0.7,3);
    \draw [red](-4,1.3) .. controls (-3.55,1.3) and (-3.2,0.95) .. (-3.2,0.5) .. controls (-3.2,0.05) and (-3.55,-0.3) .. (-4,-0.3);
    \draw [red](-2.3,-2) .. controls (-2.3,-1.55) and (-1.95,-1.2) .. (-1.5,-1.2) .. controls (-1.05,-1.2) and (-0.7,-1.55) .. (-0.7,-2);
    \draw [red](1,-0.3) .. controls (0.55,-0.3) and (0.2,0.05) .. (0.2,0.5) .. controls (0.2,0.95) and (0.55,1.3) .. (1,1.3);
    \draw [red](-4,2) .. controls (-3.45,2) and (0,-1.45) .. (0,-2);
    \draw [red](-3,3) .. controls (-3,2.45) and (0.45,-1) .. (1,-1);
    \draw [red](-4,2.5) .. controls (-3.55,2.5) and (0.5,-1.55) .. (0.5,-2);
    \draw [red](-3.5,3) .. controls (-3.5,2.55) and (0.55,-1.5) .. (1,-1.5);
    \draw [red](-4,1.75) .. controls (-3.4,1.75) and (-0.25,-1.4) .. (-0.25,-2);
    \draw [red](-2.75,3) .. controls (-2.75,2.4) and (0.4,-0.75) .. (1,-0.75);
    \draw [red](-4,2.25) .. controls (-3.5,2.25) and (0.25,-1.5) .. (0.25,-2);
    \draw [red](-3.25,3) .. controls (-3.2,2.5) and (0.5,-1.25) .. (1,-1.25);
    \draw [red](-4,3) -- (1,-2);
    \draw [mygreen, very thick](-1.5,1.7) node[fill, circle, inner sep=1.5pt] {} -- (-1.5,-0.7) node[fill, circle, inner sep=1.5pt] {};
    \node [mygreen] at (-1.25,1.5) {$\alpha$};
    
    \begin{scope}[xshift=5cm, xscale=-1]
    \draw [blue](-1.5,3) coordinate (v1) -- (-4,0.5) node (v2) {} -- (-1.5,-2) node (v3) {} -- (1,0.5) node (v4) {} -- (v1);
    \draw[blue] (v1) -- (v3);
    \draw[fill=white]  (v1) ellipse (0.5 and 0.5);
    \draw[fill=white]  (v2) ellipse (0.5 and 0.5);
    \draw[fill=white]  (v3) ellipse (0.5 and 0.5);
    \draw[fill=white]  (v4) ellipse (0.5 and 0.5);
    \draw [red](-2.5,3) .. controls (-2.5,2.5) and (-2,2) .. (-1.5,2) .. controls (-1,2) and (-0.5,2.5) .. (-0.5,3);
    \draw [red](-4,1.5) .. controls (-3.5,1.5) and (-3,1) .. (-3,0.5) .. controls (-3,0) and (-3.5,-0.5) .. (-4,-0.5);
    \draw [red](-2.5,-2) .. controls (-2.5,-1.5) and (-2,-1) .. (-1.5,-1) .. controls (-1,-1) and (-0.5,-1.5) .. (-0.5,-2);
    \draw [red](1,-0.5) .. controls (0.5,-0.5) and (0,0) .. (0,0.5) .. controls (0,1) and (0.5,1.5) .. (1,1.5);
    \draw [red](-2.05,3) .. controls (-2.05,2.7) and (-1.8,2.45) .. (-1.5,2.45) .. controls (-1.2,2.45) and (-0.95,2.7) .. (-0.95,3);
    \draw [red](-4,1.05) .. controls (-3.7,1.05) and (-3.45,0.8) .. (-3.45,0.5) .. controls (-3.45,0.2) and (-3.7,-0.05) .. (-4,-0.05);
    \draw [red](-2.05,-2) .. controls (-2.05,-1.7) and (-1.8,-1.45) .. (-1.5,-1.45) .. controls (-1.2,-1.45) and (-0.95,-1.7) .. (-0.95,-2);
    \draw [red](1,-0.05) .. controls (0.7,-0.05) and (0.45,0.2) .. (0.45,0.5) .. controls (0.45,0.8) and (0.7,1.05) .. (1,1.05);
    \draw [red](-2.15,3) .. controls (-2.15,2.65) and (-1.85,2.35) .. (-1.5,2.35) .. controls (-1.15,2.35) and (-0.85,2.65) .. (-0.85,3);
    \draw [red](-4,1.15) .. controls (-3.65,1.15) and (-3.35,0.85) .. (-3.35,0.5) .. controls (-3.35,0.15) and (-3.65,-0.15) .. (-4,-0.15);
    \draw [red](-2.15,-2) .. controls (-2.15,-1.65) and (-1.85,-1.35) .. (-1.5,-1.35) .. controls (-1.15,-1.35) and (-0.85,-1.65) .. (-0.85,-2);
    \draw [red](1,-0.15) .. controls (0.65,-0.15) and (0.35,0.15) .. (0.35,0.5) .. controls (0.35,0.85) and (0.65,1.15) .. (1,1.15);
    \draw [red](-2.3,3) .. controls (-2.3,2.55) and (-1.95,2.2) .. (-1.5,2.2) .. controls (-1.05,2.2) and (-0.7,2.55) .. (-0.7,3);
    \draw [red](-4,1.3) .. controls (-3.55,1.3) and (-3.2,0.95) .. (-3.2,0.5) .. controls (-3.2,0.05) and (-3.55,-0.3) .. (-4,-0.3);
    \draw [red](-2.3,-2) .. controls (-2.3,-1.55) and (-1.95,-1.2) .. (-1.5,-1.2) .. controls (-1.05,-1.2) and (-0.7,-1.55) .. (-0.7,-2);
    \draw [red](1,-0.3) .. controls (0.55,-0.3) and (0.2,0.05) .. (0.2,0.5) .. controls (0.2,0.95) and (0.55,1.3) .. (1,1.3);
    \draw [red](-4,2) .. controls (-3.45,2) and (0,-1.45) .. (0,-2);
    \draw [red](-3,3) .. controls (-3,2.45) and (0.45,-1) .. (1,-1);
    \draw [red](-4,2.5) .. controls (-3.55,2.5) and (0.5,-1.55) .. (0.5,-2);
    \draw [red](-3.5,3) .. controls (-3.5,2.55) and (0.55,-1.5) .. (1,-1.5);
    \draw [red](-4,1.75) .. controls (-3.4,1.75) and (-0.25,-1.4) .. (-0.25,-2);
    \draw [red](-2.75,3) .. controls (-2.75,2.4) and (0.4,-0.75) .. (1,-0.75);
    \draw [red](-4,2.25) .. controls (-3.5,2.25) and (0.25,-1.5) .. (0.25,-2);
    \draw [red](-3.25,3) .. controls (-3.2,2.5) and (0.5,-1.25) .. (1,-1.25);
    \draw [red](-4,3) -- (1,-2);
    \draw [mygreen, very thick](-1.5,1.7) node[fill, circle, inner sep=1.5pt] {} -- (-1.5,-0.7) node[fill, circle, inner sep=1.5pt] {};
    \node [mygreen] at (-1.25,1.5) {$\alpha$};
    \end{scope}
    \end{tikzpicture}
    \caption{Left: type (i), right: type (ii).}
    \label{fig:lam_in_quad}
\end{figure}

Then we get a map $\mathbf{x}^{(\tri,\ell)}:=(x_i^{(\tri,\ell)})_{i \in I}: \widehat{\ML}([\Sigma^h,f]) \to \bR^I$. 
\label{no_cusp}

For general $h \in \mathrm{Hyp}(\Sigma)$ and given $(G,\mu,\eta) \in \widehat{\ML}([\Sigma^h,f])$, we modify the surface $\Sigma^h$ by removing a small disk around each cusp and make $G$ spiralling into the direction specified by the relaxed signature $\eta$. 
 
\begin{thm}[\cite{FG07,BB09}]\label{shear1}
    For any labeled triangulation $(\tri,\ell)$ of $\Sigma$, the map
    \[
        \mathbf{x}^{(\tri,\ell)}: \widehat{\ML}(\Sigma^h)\to \bR^I,\quad (G,\mu)\mapsto (x_i^{(\tri,\ell)}(G,\mu))_{i\in I}
    \]
    is bijective.
\end{thm}
The coordinate systems $\mathbf{x}^{(\tri,\ell)}$ are called the \emph{shear coordinate systems}.
Moreover, it is known that the coordinate transformation $\mathbf{x}^{(\tri',\ell)} \circ (\mathbf{x}^{(\tri,\ell)})^{-1}$ associated with a flip coincides with the tropical $\X$-transformation \cite{FG07}. Indeed, it follows from a computation on the rational part as studied in \cite{FG07} and a continuity argument. Therefore, the pairs $(\ve^{(\tri,\ell)},\mathbf{x}^{(\tri,\ell)})$ associated with labeled triangulations $(\tri,\ell)$ can be view as tropical seeds in the semifield $\bR^\trop$. The shear coordinate systems combine to give  a canonical isomorphism $\widehat{\ML}(\Sigma^h) \cong \X_\Sigma(\bR^\trop)$ that is compatible with \eqref{eq:graph_emb}. 

\smallskip
\paragraph{\textbf{A trivialization of the lamination bundle.}} Bunching the $MC(\Sigma)$-equivariant isomorphisms $\widehat{\ML}(\Sigma^h) \cong \X_\Sigma(\bR^\trop)$ for all $h \in \mathrm{Hyp}(\Sigma)$, we get an $MC(\Sigma)$-equivariant trivialization
\begin{align*}
    \eML_\Sigma \cong \widehat{\cT}(\Sigma) \times \X_\Sigma(\bR^\trop)
\end{align*}
of topological fiber bundles over $\widehat{\cT}(\Sigma)$. The isomorphism $\widehat{\cT}(\Sigma) \cong \X_\Sigma(\pos)$ further allows us to obtain an $MC(\Sigma)$-equivariant homeomorphism
\begin{align}\label{eq:bundle_trivial}
    \eML_\Sigma \cong \X_\Sigma(\pos) \times \X_\Sigma(\bR^\trop).
\end{align}
It is also equivariant under the action of $\bR_{>0}$ rescaling the measured geodesic laminations and the shear coordinates.

\section{A brief review on the earthquake maps}\label{sec:eq}
We assume that $\Sigma = (\Sigma,M)$ is a marked surface as in \cref{cl_surf}. 

\subsection{Earthquake theorem}
We will construct 
the earthquake map
\[
    E\cl \eML_\Sigma \to \widehat{\cT}(\Sigma),
\]
following \cite{BB09,BKS16,EM,A18}. 
We use 
the description of $\widehat{\cT}(\Sigma)$ given in \cref{prop:enhanced_Farey}.
Let $([\rho^h,\psi^h],[G,\mu,\eta])\in\eML_\Sigma$, and $\Sigma^h=\HHH/\rho^h(\pi_1(\Sigma))$ the underlying hyperbolic surface. Let $\widetilde{\Sigma^h} \subset \HHH$ be the universal cover of $\Sigma^h$, and $(\widetilde{G}, \widetilde{\mu})$ be the full preimage of $(G,\mu)$.

We will describe the earthquake map $\widetilde{E}_{(\widetilde{G},\widetilde{\mu})}$ in the universal cover. We call each connected component of $\widetilde{\Sigma^h}\setminus\widetilde{G}$ a \emph{gap} of $\widetilde{G}$. We fix a gap $S_0$.
We will endow each gap $S$ with a hyperbolic element $E_S\in\PSL$. First, we endow $S_0$ with the identity. For each gap $S$, we take an arc $\alpha_S$ transverse to $\widetilde{G}$, from $S_0$ to $S$ and never revisiting the same gap. We take $m\in\bZ_{\geq 1}$ and a partition $x_0,...,x_k$ of $\alpha_S$ satisfying the following conditions:

\begin{itemize}
    \item The point $x_0$ is the endpoint of $\alpha_S$ in $S_0$. The point $x_k$ is the endpoint of $\alpha_S$ in $S$. 
    \item Each point $x_i$ is in a gap ($i=1,...,k-1$),
    \item Let $[x_i,x_{i+1}]$ be the segment from $x_i$ to $x_{i+1}$ in $\alpha_S$. The hyperbolic length of $[x_i,x_{i+1}]$ is less than $1/m$ $(i=0,...,k-1)$.    
\end{itemize}

If there are  leaves of $\widetilde{G}$ transverse to $[x_i,x_{i+1}]$, we take a leaf  $l_i$ from among them and define $\widetilde{E_{i}}$ as the hyperbolic element whose axis is $l_i$, whose translation distance is $\widetilde{\mu} ([x_i,x_{i+1}])$ ($i=0,...,k-1$) and whose direction is to the left (i.e. the restriction $\widetilde{E_{i}}|_{S}$ of $\widetilde{E_{i}}$ moves $S$ to the left from $S_0$'s perspective).  Otherwise, we define $\widetilde{E_{i}}$ as the identity. We define $\widetilde{E_{S,m}}:=\widetilde{E_{k-1}}\circ \widetilde{E_{k-2}}\circ\cdots\circ \widetilde{E_{0}}$. When $m\to \infty$, the hyperbolic element $\widetilde{E_{S,m}}$ converges in $\PSL$ regardless of how $\alpha_S$, $x_0,...,x_k$ and $l_0,...,l_k$ are chosen (\cite[Lemma II.3.4.4]{EM}). We endow any gap $S$ with the limit $\dsp \widetilde{E_S}:=\lim_{m\to \infty}\widetilde{E_{S,m}}$.

We set $\widetilde{E}_{(\widetilde{G},\widetilde{\mu})}:=\dsp\bigcup_{S\text{:\,gap}}\widetilde{E_S}:\widetilde{\Sigma^h}\setminus \widetilde{G}\to \HHH$.
It is continuously extended to a map $\partial_\infty\widetilde{E}_{(\widetilde{G},\widetilde{\mu})}:\partial_\infty\widetilde{\Sigma^h}\to\SSS$, where $\partial_\infty\widetilde{\Sigma^h}$ is the intersection of the closure of $\widetilde{\Sigma^h}$ and $\SSS$.

\begin{dfn}
    We define the \emph{earthquake map}
    \begin{align}\label{eq:EQ}
        E: \eML_\Sigma \to \widehat{\cT}(\Sigma), \quad E([\rho^h,\psi^h],[G,\mu,\eta])=[\rho^{h'},\psi^{h'}],
    \end{align}
    as follows. 
    \begin{enumerate}
        \item We define a homomorphism $\rho^{h'}:\pi_1(\Sigma)\to \PSL$ by the condition $\partial_\infty \widetilde{E}_{(\widetilde{G},\widetilde{\mu}))}\circ\rho^h(\gamma)=\rho^{h'}(\gamma)\circ\partial_\infty \widetilde{E}_{(\widetilde{G},\widetilde{\mu}))}$ on $\partial_\infty\widetilde{\Sigma^h}$ for any $\gamma\in\pi_1(\Sigma)$. 
        Then $\rho^{h'}$ is a Fuchsian representation by \cite[Proposition 3.2.]{BKS16}. 
        \item Define $\psi^{h'}:\cF_\infty(\Sigma) \to \SSS$ by $\psi^{h'}:=\partial_\infty \widetilde{E}_{(G,\mu)}\circ \psi^h$ if $\eta_p=+$ for all $p \in M_\circ$. Then $\psi^{h'}$ is $\rho^{h'}$-equivariant and order-preserving. Otherwise, we need to modify $\psi^{h'}$. See \cite{BB09,BKS16}.
        
\end{enumerate}
\end{dfn}

Via the homeomorphism \eqref{eq:bundle_trivial}, we can rewrite \eqref{eq:EQ} as
\begin{align}\label{eq:EQ_split}
    E\cl  \cX_\Sigma (\bR_{>0})\times\cX_\Sigma (\bR^{\trop})\rightarrow\cX_\Sigma (\bR_{>0}).
\end{align}
In the latter setting, the earthquake theorem is stated as follows:

\begin{thm}[Earthquake theorem {\cite{Thu84,BKS16}}]
    
    If $M_\partial=\emptyset$, the map 
    \[
        E(g_0,-)\cl\cX_\Sigma (\bR^{\trop})\to \cX_\Sigma (\pos)
    \]
    is a 
    homeomorphism for any $g_0 \in \X_\Sigma(\pos)$.
\end{thm}
We believe that this theorem holds for any marked surface. 
When we fix a measured geodesic lamination $L \in \X_\Sigma(\bR^\trop)$ instead, the automorphism $E(-,L)$ on $\widehat{\cT}(\Sigma)=\X_\Sigma(\pos)$ is called the \emph{earthquake
} along $L$.

\subsection{Explicit formula for earthquakes along ideal arcs}\label{subsec:eq_exp}

Let $\gamma$ be an ideal arc on $\Sigma$ connecting two marked points $p_1,p_2$. Let $L(\gamma) \in \X_\Sigma(\bR^\trop)$ be the point characterized by the condition $x_i^{(\tri,\ell)}(L(\gamma))=\delta_{ik}$ for $i \in I$ and any labeled triangulation $(\tri,\ell)$ such that $\ell(k)=\gamma$. Notice that this condition does not depend on the choice of $(\tri,\ell)$. For $[\Sigma^h,f] \in \widehat{\cT}(\Sigma)$, the corresponding measured geodesic lamination $L^h(\gamma) \in \eML([\Sigma^h,f])$ is described as follows:
\begin{itemize}
    \item the support is given by $f(\gamma)$;
    \item the transverse measure is given by $\# (f(\gamma) \cap \alpha)$ 
    for any transverse arc $\alpha$;
    \item the relaxed signature is given by $\eta_p:=1$ if $p \in \{p_1,p_2\}$, and $\eta_p:=0$ otherwise. 
\end{itemize}

\begin{prop}\label{eq_surf}
    Let $\gamma$ be an ideal arc, and $t >0$ a positive number. 
    Then for any $g=[\Sigma^h,f] \in\cT(\Sigma)$, we have
    \[
        X^{(\tri,\ell)}_i(E(g,t\cdot L^h(\gamma)))= e^{t \delta_{ik}}\cdot X^{(\tri,\ell)}_i(g)
    \]
    for $i \in I$. Here $(\tri,\ell)$ is any labeled triangulation $(\tri,\ell)$ such that $\ell(k)=\gamma$.  
\end{prop}
\begin{proof}
    For $i \in I$, let $\square_i$ denote the unique quadrilateral in $\tri$ that contains $\ell(i)$ as a diagonal. For a point $[\Sigma^h,f] \in \widehat{\cT}(\Sigma)$, choose a geodesic lift $\widetilde{\tri} \subset \widetilde{\Sigma^h} \subset \HHH$ (tesselation) of $f(\tri)$ with respect to the universal covering $\widetilde{\Sigma^h} \to \Sigma^h$ determined by the hyperbolic structure $h$. We may choose this lift so that the ideal vertices of a lift $\widetilde{\square_k}$ of $\square_k$ are $0,1,\infty,-X_k$ in this counter-clockwise order and $\{0,\infty\}$ are the ideal points of $\widetilde{f(\gamma)}$, where $X_k:=X_k^{(\tri,\ell)}(g)$ is the cross-ratio coordinate of $g$ for $k \in I$. 
    Moreover, from the construction, we may normalize the earthquake $E_{t\cdot L^h(\gamma)}$ so that its restriction to the triangle $(0,1,\infty)$ is identity. Then its restriction to the other triangle $(\infty,-X_k,0)$ is the hyperbolic element in $\PSL$ with translation length $t$ and axis $\widetilde{f(\gamma)}=\sqrt{-1} \cdot \bR_{>0}$. Explicitly, it is given by $\begin{pmatrix} e^{t/2} & 0 \\ 0 & e^{-t/2} \end{pmatrix}$. 
    Then we get
    \[
        X_k^{(\tri,\ell)}(E(g,t\cdot L^h(\gamma)))=[0:1:\infty:-e^t\cdot X_k]=e^t\cdot X_k^{(\tri,\ell)}(g).
    \]
    For $i \neq k$, notice that a common hyperbolic element is applied entirely on any lift of $\square_i$. Hence $X_i^{(\tri,\ell)}(E(g,t\cdot L^h(\gamma)))= X_i^{(\tri,\ell)}(g)$, since the cross-ratio is an $PSL_2(\bR)$-invariant of four points. Thus the assertion is proved.
\end{proof}

Since the earthquake deformations along two measured geodesic laminations with disjoint supports commute with each other, we obtain the following:

\begin{cor}
    Let $L \in \cC_{(\tri,\ell)}^+$ for some labeled triangulation $(\tri,\ell)$ (\cref{def:positive_cone}), and $t_i:=x_i^{(\tri,\ell)}(L) \geq 0$ its shear coordinates for $i \in I$. For $g=[\Sigma^h,f]\in \widehat{\cT}(\Sigma)$, let $L^h \in \eML([\Sigma^h,f])$ be the corresponding point. Then we have
    
    \[
        X^{(\tri,\ell)}_i(E(g,L^h))= e^{t_i}\cdot X^{(\tri,\ell)}_i(g)
    \]
    for $i \in I$. 
\end{cor}
This corollary gives a coordinate description of the earthquake map on $\cX_\Sigma(\bR_{>0})\times\cC_{(\tri,\ell)}^+$.


\end{document}